\documentclass[12pt]{amsart}
\usepackage{amsmath}
\usepackage{amssymb}
\usepackage{amscd}
\usepackage{epsfig}
\usepackage{amsthm}

\def\t{\tilde}

\newcommand{\length}{\operatorname{\ul{length}_{\Delta}}}
\newcommand{\Length}{\operatorname{length_{\Delta}}}

\newtheorem{dfn}{Definition}[section]
\newtheorem{rem}[dfn]{Remark}
\newtheorem{thm}[dfn]{Theorem}
\newtheorem{defn}[dfn]{Definition}
\newtheorem{lem}[dfn]{Lemma}
\newtheorem{lemma}[dfn]{Lemma}

\newtheorem{prop}[dfn]{Proposition}

\newtheorem{cor}[dfn]{Corollary}
\newtheorem{conj}[dfn]{Conjecture}
\newtheorem{conjecture}[dfn]{Conjecture}

\newtheorem{notation}[dfn]{Notation}
\newtheorem{ex}[dfn]{Example}

\newtheorem{convention}[dfn]{Convention}
\newtheorem{observation}[dfn]{Observation}

\def\proof{\par\medskip\noindent{\it Proof: }}

\def\CR{\curvearrowright}
\def\acts{\CR}

\def\embed{\hookrightarrow}

\def\P{{\mathcal P}}
\def\C{{\mathbb C}}
\def\R{{\mathbb R}}

\def\Z{{\mathbb Z}}
\def\K{{\mathbb K}}

\def\E{{\mathcal E}}

\def\O{{\mathcal O}}

\def\Q{{\mathbb Q}}

\def\N{{\mathbb N}}

\def\eps{\epsilon}
\def\al{\alpha}
\def\be{\beta}
\def\ga{\gamma}

\def\De{\Delta}

\def\Del{\Delta}
\def\Si{\Sigma}
\def\si{\sigma}

\def\la{\lambda}
\def\La{\Lambda}

\def\om{\omega}

\def\>{\rangle}
\def\<{\langle}
\def\del{\delta}
\def\D{\partial}
\def\3{\ss}
\def\8{\infty}

\def\ol{\overline}
\def\ul{\underline}
\def\ov{\overrightarrow}
\newcommand{\restr}{\mbox{\Large \(|\)\normalsize}}

\title{Structure of the tensor product semigroup}

\begin{document}

\title{Structure of the tensor product semigroup}
\author{Misha Kapovich and John J. Millson}
\dedicatory{To the memory of S. S. Chern}
\date{August 10, 2005}

\begin{abstract}
We study the structure of semigroup $Tens(G)$ consisting of triples of dominant weights $(\la,\mu,\nu)$ 
of a complex reductive Lie group $G$ such that
$$
(V_\la\otimes V_\mu
\otimes V_\nu)^G \ne 0. 
$$ 
We prove two general structural results for $Tens(G)$ and give an explicit computation of $Tens(G)$ 
for $G=Sp(4,\C)$ and $G=G_2$.  
\end{abstract}

\maketitle

\section{Introduction}

Suppose that $G$ is a complex reductive Lie group.
Finite-dimensional irreducible representations $V_\la$ of $G$ are
parameterized by their highest weights $\la\in \Del\cap L$, where
$\De$ is the positive Weyl chamber, $L$ is the character lattice
of a maximal (split) torus in $G$. One of the basic questions of
the representation theory is to decompose tensor products
$V_\la\otimes V_\mu$ into sums of irreducible representations.
Accordingly, we define the set
$$
Tens(G):= \{(\la,\mu,\nu)\in (\Del \cap L)^3 : (V_\la\otimes V_\mu
\otimes V_\nu)^G \ne 0\}.
$$
For simply-connected Lie groups with the root system $R$ we will
write $Tens(R)$ instead of $Tens(G)$. It was known for a long time,
see for example \cite[Theorem 9.8]{KLM3}, 
that the set $Tens(G)$ forms a semigroup with respect to the
addition. The goal of this paper is to provide more specific
structural theorems for $Tens(G)$ and to make an explicit
computation of $Tens(Sp(4,\C))$ and $Tens(G_2)$.

\begin{thm}
For each complex reductive Lie group $G$, the set $Tens(G)$ is a finite union of
elementary subsets of $L^3$.
\end{thm}

Here an {\em elementary subset} is a subset given by a finite
system of linear inequalities (with integer coefficients) and
congruences. Thus, to ``describe'' $Tens(G)$ one would have to
find these inequalities and congruences. The above theorem is an
analogue of a theorem by C.~Laskowski \cite{Laskowski}, who proved
a similar statement for the structure constants of spherical Hecke
rings. 

\medskip
Our next theorem provides a glimpse of what these inequalities and
congruences might look like. In \cite{BS} and \cite{KLM1} there have been defined a
finite-sided homogeneous polyhedral cone $\P(G)=D_3(G/K)\subset \Del^3$
(where $K$ is a maximal compact subgroup of $G$), given
by the inequalities of the form: 
$$
\<\varpi_i, w_1\la\> + \<\varpi_i, w_2 \mu\> +\<\varpi_i,
w_3 \nu\> \le 0,
$$
where $\varpi_i$ are the fundamental weights of $R$ and $w_i$ are certain elements of 
the Weyl group of $G$ determined  by the ``Schubert calculus''.
It is known (see for instance \cite{KLM3}) that
$$
Tens(G)\subset \P(G) \cap \{ \si=(\la, \mu, \nu)\in L^3:
Tr(\si)\in Q(R)\}. 
$$
Here and in what follows $Q(R)$ is the root lattice and $Tr(\si):= \la+\mu+\nu$.

\begin{thm}
\label{thm:deep}
There exists a vector $\si\in \P(G)$ such that
$$
(\P(G)+\si) \cap Tens(G)= (\P(G)+\si) \cap  \{ (\la, \mu, \nu)\in L^3: Tr(\si)\in
Q(R)\}.
$$
\end{thm}

In other words, inside of the translated cone $\P(G)+\si$ the {\em
necessary} conditions on $(\la,\mu,\nu)$ to belong to $Tens(G)$ are
also {\em sufficient}. Observe that the intersection $\{ \si\in L^3: Tr(\si)\in Q(R)\}\cap (\P(G)+\si)$ is an
elementary subset of the lattice $L^3$. The length of the vector
$\si$ can be explicitly estimated. In section \ref{subcone} we
make such an estimate in the case when $R=B_\ell$; the same methods 
work for other root systems.  

\medskip
In section \ref{B2C2} we will compute the semigroup $Tens=Tens(G)$
for the group $G=Sp(4,\C)$ and show that this set is {\em not}
an elementary set itself. Below
%$K=SU(2)$ is a maximal compact
%subgroup in $G$,
$P(R)=L$ is the weight lattice of $G$ and
$\varpi_2$ is the longest fundamental weight. We let
$$
\La:= \{\si\in L^3: Tr(\si)\in Q(R)\}.
$$
%We first write the description of $Tens$ in English and then
%translate this description into decomposition of $Tens$ into
%elementary subsets:

Then 

\begin{thm}
$\si=(\la,\mu,\nu)\in \P(G)\cap \La$ belongs to $Tens$ if and only
if one of the following is satisfied:

1. At least two of the vectors $\la, \mu, \nu$ are not multiples
of $\varpi_2$.

2. If two of the vectors $\la, \mu, \nu$ are multiples of
$\varpi_2$ then $Tr(\si) \in 2 P(R)$.
\end{thm}

In section \ref{sect:BC} we express $Tens$ as a union of 6 elementary subsets. 
We refer the reader to \cite{KLM1, KuLM} (see also section \ref{inequalities} of this paper) 
for the explicit computation of the polyhedra $\P(G)$ for all 
complex simple Lie group of rank $\le 3$. 

In section \ref{G2} we present a computation of the semigroup $Tens(R)$ for the root system $R=G_2$.  

\medskip
We note that prior to the present paper, the complete description of the semigroup $Tens(G)$ was 
known only for the groups of type $A_n$.  In this case the Saturation Theorem 
of Knutson and Tao \cite{KT} implies that  
$$
Tens(G)= \P(G) \cap \{ \si\in L^3:
Tr(\si)\in Q(R)\}.
$$

Below is a conjecture on the structure of the semigroup $Tens(G)$ for an arbitrary root system $R$. Here
$L$ denotes the character lattice of a maximal torus. A triple of dominant weights 
$\si=(\la_1,\la_2,\la_3)$ is called {\em singular} if at least one of the vectors $\la_i$ is singular, i.e. 
belongs to a wall of the Weyl chamber $\De$. 

\begin{conjecture}
\label{C1}
1. Suppose that $R$ is simply-laced. Then
$$
Tens(G)=\{\si=(\la_1,\la_2,\la_3)\in \P(G): \la_i\in L, Tr(\si)\in Q(R)\}.
$$

2. Suppose that $\si=(\la_1,\la_2,\la_3)\in \P(G)\cap L^3$ is a nonsingular triple. Then
$\si\in Tens(G)$ if and only if $\si\in \P(G)$ and $Tr(\si)\in Q(R)$.

3. Suppose that $P(R)=Q(R)$. Then in the decomposition of $Tens(G)$ as the union of elementary sets,
the elementary sets are given by
inequalities only and there are no congruence conditions.
\end{conjecture}

The above conjecture holds for the root systems $A_n$, $B_2=C_2$ and $G_2$, it is also supported by a number of 
computer-based calculations with the higher rank root systems. 

\begin{rem}
Let $T^{reg}$ be the set of nonsingular triples $\si\in \P(G)\cap L^3$.
It was observed by J.~Bernstein that in the decomposition of $T^{reg}$ into elementary sets, there are no congruence
conditions apart from the ``obvious'' condition $Tr(\si)\in Q(R)$.
\end{rem}

A less ambitious conjecture (which actually follows from either Part 1 or 2 of Conjecture \ref{C1}) is

\begin{conjecture}
[S. Kumar]
\label{Kumar2}
There exists a triple $\si\in \P(G)\cap L^3, Tr(\si)\in Q(R)$ such that $\si\notin Tens(G)$ if and only if
there exists a singular triple $\si$ with the above properties.
\end{conjecture}

\medskip
{\bf Acknowledgements.} During the work on this paper the first
author was supported by the NSF grant DMS-04-05180; part of this work was done when 
he was visiting the Max Plank Institute for Mathematics in Bonn. The second
author was supported by the NSF grant DMS-04-05606. The authors
gratefully acknowledge support of these institutions. 

The possibility of Theorem \ref{thm:deep} was suggested by J.~Bernstein
in a conversation at Oberwolfach. T.~Haines and S.~Kumar have told us about about Theorem \ref{kumar}.
The authors are grateful to J.~Bernstein, T.~Haines, S.~Kumar and C.~Laskowski for 
these and other useful conversations. 
The second author would like to acknowledge how exciting it was to
learn Chern-Weil theory and the theory of Chern-Simons invariants first-hand
when he was a graduate student in Berkeley in the early seventies.

\section{Review of the path model for the representation theory of complex
reductive Lie groups}

\subsection{Root systems and Coxeter complexes}

Let $V$ be a  finite-dimensional Euclidean vector space and $R\subset V$ be a root system in $V$.
Then the collection of coroots $R^\vee$ determines a root system in $V^*$.
Using the metric on $V$ we will be identifying $R$ and $R^\vee$ with root systems in $V$.
Thus we will think of both $R$ and $R^\vee$ as
linear functionals on $V$. Given $R$ we define the {\em affine Coxeter group} $W_{aff}:=W_{R^\vee,aff}$
as the group generated by reflections in the {\em walls}
$$
H=\{x\in V: \al(x)=t\}, \al\in R^\vee, t\in \Z.
$$

We let $W=W_{R^\vee,sph}=W_{R,sph}$ denote the linear part of $W_{R,aff}$,
which is the same as the stabilizer of the origin in $W_{R^\vee,aff}$.
The group $W$ is the Weyl group
of the root system $R$ (and $R^\vee$).

\begin{rem}
In the context of the representation theory of a complex semisimple
Lie group $G$, the space $V$ equals $X^*(T)\otimes \R$, where $T$ is a
maximal torus in $G$.
Thus $R\subset V$ and the walls in $V$ are given by the coroot system $R^\vee$.
\end{rem}

A {\em Weyl chamber $\De$} is a fundamental domain
for the finite reflection group $W$, it is bounded by certain walls passing through the
origin.  Given a vector $v\in \De$ we let $v^*\in \De$ denote the {\em contragredient vector}
$v^*=w_0(-v)$, where $w_0\in W$ is the longest element, i.e. the element which
carries $-\De$ to $\De$.

The group $W_{aff}$ acts by isometries on the Euclidean space $V$. Let $A$ denote
the affine space
underlying $V$. The pair $(A,W_{aff})$ is called a {\em Euclidean
Coxeter complex}. Let ${\mathcal W}$ denote the union of all
walls. Then the closures of connected components  of $A\setminus
{\mathcal W}$ are called {\em alcoves}. If $R$ is irreducible and spans $V$
then alcoves are simplicies and therefore  $(A,W_{aff})$ has
natural structure of a simplicial complex.

Each alcove is a fundamental domain for the action $W_{aff}\acts A$. Pick an alcove 
$a$ and call it a fundamental alcove. We have a natural projection $\theta: A\to a$ sending each 
point $v\in A$ to the inique intersection point $\theta(v)\in W_{aff}\cdot v\cap a$. The image $\theta(v)$ 
is called the {\em type} of $v$.

{\em Special verticies} of the complex $(A,W_{aff})$ are points whose stabilizer in
$W_{aff}$ is isomorphic to  $W_{sph}$. The weight group $P(R)$ acts
simply transitively on the set of special verticies.

Let $h: A\to A$ be a dilation of $A$, i.e. an affine map of the
form $x\mapsto kx+b$, where $b\in V$, $k>0$. The number $k$ is the
{\em conformal factor} of $h$. We define $Dil(A, W_{aff})$ to be the
semigroup of dilations $h$ of $A$ such that
$$
hW_{aff}h^{-1}\subset W_{aff}.
$$
Then each $h\in Dil(A,W_{aff})$ sends verticies of $(A,W_{aff})$ to
verticies of $(A,W_{aff})$, walls to walls, etc.

\subsection{Chains}

Let $R\subset V$ be a root system in a Euclidean vector space $V$,
$W$ be the Weyl group of $R$. We pick a Weyl chamber $\Del$ for
$W$, this determines the set of positive roots $R_+$ and the set of the simple roots $\Phi$ in $R$,
as well as positive and simple coroots. The following notion of {\em chains} and the {\em partial order
$\ge$ on $V\setminus \{0\}$} was introduced by P.~Littelmann in \cite{Littelmann2}.

\begin{defn}\label{D0}
 A $W$--{\em chain} in $V$ is a finite sequence $(\eta_0,...,\eta_m)$ of nonzero
vectors in $V$ so that for each $i=1,...,m$ there exists a positive
coroot $\be_i\in R^\vee$ so that the corresponding reflection
$\tau_i:=\tau_{\be_i}\in W$ satisfies

1. $\tau_i(\eta_{i-1})= \eta_{i}$.

2. $\be_i(\eta_{i-1}) <0$.
\end{defn}

Then $\eta\ge_W \nu$ if there exists a chain from $\eta$ to $\nu$.
Most of the time we will abbreviate $\ge_W$ to $\ge$. We say that a chain
$$
\eta_0'\ge ...\eta_i' \ge ... \ge \eta'_n$$
 is a {\em refinement} of the chain
$$
\eta_0\ge ... \ge \eta_i \ge ... \ge \eta_k
$$
if $\eta_0=\eta_0', \eta_n'=\eta_k$ and
$$
\{\eta_0, ...,\eta_k\}\subset \{\eta_0', ...,\eta'_n\}.
$$
A chain which does not admit a proper refinement is called {\em maximal}.

Define a (nontransitive) relation $\sim=\sim_W$ on
$V\setminus \{0\}$ by
$$
\mu \sim_W \nu \iff
$$
\centerline{$\mu, \nu$ belong to the same Weyl chamber of $W$.}

We  write $\la\gtrsim \xi$ if there exist $\mu, \nu$ so that
$$
\la
\ge \mu \sim \nu \ge  \xi.$$
 Accordingly, we define {\em
generalized chains} as sequences of nonzero vectors in $V$:
$$
\eta_0\gtrsim \eta_1 \gtrsim ... \gtrsim \eta_m.
$$

\subsection{Buildings}
\label{inequalities}

Our discussion of buildings follows \cite{KleinerLeeb}. We refer the reader
to \cite{Brown}, \cite{Ronan}, \cite{Rousseau} for the more combinatorial discussion.

Fix a spherical or Euclidean (discrete) Coxeter complex $(A,W)$, where $A$ is a
Euclidean space $E$ or a unit sphere $S$ and $W=W_{aff}$ or $W=W_{sph}$
is a discrete Euclidean or a spherical Coxeter group acting on $A$.

A metric space $Z$ is called {\em geodesic} if every pair of points $x, y$ in $Z$ can be connected by a
geodesic segment $\ol{xy}$.

Let $Z$ be a metric space. A {\em geometric structure} on $Z$ {\em
modeled on} $(A,W)$ consists of an atlas of isometric embeddings
$\varphi:A\embed Z$ satisfying the following compatibility
condition: For any two charts $\varphi_1$ and $\varphi_2$, the
transition map $\varphi_2^{-1}\circ\varphi_1$ is the restriction of
an isometry in $W$. The charts and their images,
$\varphi(A)=a\subset Z$, are called {\em apartments}. We will
sometimes refer to $A$ as the {\em model apartment}. We will require
that there are {\em plenty of apartments} in the sense that any two
points in $Z$ lie in a common apartment. All $W$-invariant notions
introduced for the Coxeter complex $(A, W)$, such as rank, walls, singular
subspaces, chambers etc., carry over to geometries modeled on $(A,
W)$. If $a, a'\subset X$ are alcoves (in the Euclidean case) or chambers
(in the spherical case) then there exists an apartment $A'\subset X$ containing $a\cup a'$:
Just take regular points $x\in a, x'\in a'$ and an apartment $A'$ passing through $x$ and $x'$.

A geodesic metric space $Z$ is said to be a $CAT(0)$-space (resp, $CAT(1)$-space) if
{\em geodesic triangles in $Z$ are ``thinner'' than geodesic
triangles in $\R^2$} (resp. in the unit sphere $S^2$). We refer the reader to \cite{Ballmann}  for the
precise definition.

\begin{defn}
A Euclidean (resp. spherical) building is a $CAT(0)$-space (resp. $CAT(1)$-space) modeled on a
Euclidean (resp. spherical) Coxeter complex.
\end{defn}

A building is called {\em thick} if every wall is an intersection of
apartments. A non-thick building can always be equipped with a
natural structure of a thick building by reducing the Coxeter group.

\medskip
Let $\K$ be a local field with a (discrete) valuation $\nu$ and valuation ring $\O$. Given a
split reductive algebraic group $\ul{G}$ over $\Z$, and a nonarchimedian Lie
group $G=\ul{G}(\K)$ we can associate with it a Euclidean building (a
Bruhat-Tits building) $X=X_G$. We refer the reader to \cite{BT},
\cite{KLM3} and \cite{Rousseau} for more detailed discussion of the
properties of $X$. Here we only recall that:

1. $X$ is thick and locally compact.

2. $X$ is modeled on a Euclidean Coxeter complex $(A,W_{aff})$ whose
dimension equals the rank of $\ul{G}$, and the root system is
isomorphic to the root system of $\ul{G}$.

3. $X$ contains a special vertex $o$ whose stabilizer in $G$ is $\ul{G}(\O)$.

\begin{ex}
Let $X$ be a (discrete) Euclidean building, consider the  {\em
spaces of directions} $\Si_x X$. We will think of this space as the
space of germs of non-constant geodesic segments $\ol{xy}\subset X$.
As a polysimplicial complex $\Si_xX$ is just the link of the point
$x\in X$. The space of directions has the structure of a spherical
building modeled on $(S,W_{sph})$, which is thick if and only if $x$
is a special vertex of $X$, see \cite{KleinerLeeb}. The same applies
in the case when $X$ is a spherical building.
\end{ex}

\medskip
Let $(A,W_{aff})$ be a Euclidean Coxeter complex and pick a Weyl chamber $\De\subset A$.

Given a pair of points $x, y\in A$ we define their $\De$-distance $d_\De(x,y)$ by taking the
vector $v:= y-x$ and applying to it an element $w\in W_{sph}$ such that  $u:= w(v)\in \De$. Then $d_\De(x,y):=u$.

Suppose that $X$ is a Euclidean building modeled on $(A,
W_{aff})$. We define a $\De$-distance in $X$ as follows. For a
pair of points $x, y\in X$ pick an apartment $\phi: A\to A'\subset
X$ such that $A'$ contains $x, y$. Then set
$$
d_\De(x,y):= d_\De(\phi^{-1}(x),\phi^{-1}(y)).
$$
It is easy to see that this distance is independent of the choice
of $\phi$. A similar definition applies if $X=G/K$ is a nonpositively curved symmetric space,
where $A$ is a maximal flat in $X$ and $W_{sph}\acts A$ is the Weyl group of $G$.

\medskip
{\bf Generalized triangle inequalities.} Suppose that $X$ is a nonpositively curved symmetric space or a
Euclidean building as above. Define the set
$$
D_3(X):=\{(\la,\mu,\nu)\in \De^3: \exists \hbox{~a geodesic triangle~} [x,y,z]\subset X \hbox{~with~}
$$
$$
d_\De(x,y)=\la, d_\De(y,z)=\mu, d_\De(z,x)=\nu\}.
$$
It is proven in \cite{KLM1, KLM2} that $D_3(X)$ is a convex homogeneous polyhedral cone which
depends only on the pair $(A, W_{sph})$ and nothing else, therefore we will frequently use
the notation $\P(G)$ for $D_3(G/K)$, where $G$ is a reductive Lie group with a maximal compact subgroup
$K$. In the case when $G$ is a complex semisimple Lie group, the inequalities defining this polyhedron have the form

1. Stability inequalities $\psi_j(\la,\mu,\nu)\ge 0$:
$$
-\<\varpi_i, w_1\la\> - \<\varpi_i, w_2 \mu\> -\<\varpi_i,
w_3 \nu\> \ge 0,
$$
where $\varpi_i$ are the fundamental weights of $R$ and $w_i$ are certain elements of the Weyl group of $G$.

2. Chamber inequalities $\psi_k(\la,\mu,\nu)\ge 0$:
$$
\al^\vee(\la)\ge 0, \al^\vee(\mu)\ge 0, \al^\vee(\nu)\ge 0,
$$
where $\al$ are simple roots in $R$.

\medskip 
In \cite{KLM1} the polyhedra $D_3(G)$ were computed for all complex semisimple Lie groups $G$
of rank $2$. Below we provide the explicit set of {\em stability inequalities} for this polyhedron in the case 
$G=Sp(4,\C)$. 

The Weyl chamber $\Delta$ is given by
$$
\Delta = \{(x,y): x\ge y\ge 0 \}.$$ 
We will omit these inequalities from the list of inequalities defining $\P(G)$ and will list only the {\em stability
inequalities}. Instead of the notation $(\la, \mu, \nu)$ for elements of $\De^3$, we will 
use the more symmetric notation  $(\la_1, \la_2, \la_3)$, where $\la_i=(x_i, y_i), i=1, 2,3$. 

The system of stability inequalities defining $\P(Sp(4,\C))$ breaks into two subsystems (since $Sp(4,\C)$ has rank 2). 
The first subsystem  is given by
\begin{align*}
x_i \leq x_j + x_k , \quad \{i,j,k\}= \{1,2,3\} \\
y_i \leq y_j + x_k ,  \quad \{i,j,k\}= \{1,2,3\}.
\end{align*}
In order to describe the second subsystem we set 
$$
S = x_1 + y_1 + x_2 + y_2 + x_3 + y_3.$$
The second subsystem is then given by
$$
x_i  + y_j \leq S/2, \quad 1\leq i,j \leq 3. $$

\subsection{LS paths and their generalizations}

Suppose that $(A, W_{aff})$ is a Euclidean Coxeter complex.
Given a point $x\in A$ let $W_x$ denote the stabilizer of $x$ in $W_{aff}$.
For a vector $v$ in $V$ define the path $\pi_v$ by the formula:
$$
\pi_v(t)=tv, \quad t\in [0,1].
$$

In what follows we will assume that all paths are (re)parameterized to have constant speed 
and domain $I:=[0,1]$.

Given two paths $p_1, p_2$ in $A$, we define their {\em
concatenation} $p=p_1* p_2$ by
$$
p(t)= \left\{
\begin{array}{c}
p_1(t), \quad t\in [0, 1],\\ p_1(1)-p_2(0)+ p_2(t), \quad t\in [1,
2].
 \end{array}\right.
$$

\medskip
Suppose that $p: I\to A$ is a path and $J=[a,b]$ is nondegenerate
subinterval in $I=[0,1]$. We will use the notation  $p\restr_J$ to denote the 
%function-theoretic 
restriction of $p$ to $[a,b]$ 
%(we {\em do not} reparameterize this restriction 
%to have the domain equal to $I$). 

If $p$ is a PL path in $A$ which is the concatenation
$$
\ol{x_1 x_2} * ...* \ol{x_{n} x_{n+1}}
$$
of geodesic segments,
then the $\De$-length of $p$, denoted $\Length(p)$, is the sum
$$
\sum_{i=1}^n d_{\De}(x_i, x_{i+1}).
$$

Given a PL path $p$ in $A$ we use the notation $p'_{\pm}(t)\in T_{p(t)}(A)$ to denote the
derivatives of $p$ at $t$ from the left and from the right.

We say that a path $p: [0,1]\to A$ is a {\em billiard path}
if for each $t\in [0,1]$,% the vectors
$$
p'_-(t)\in W_{p(t)} p'_+(t).
$$
If $p:=p_1* ...* p_m$ is a concatenation of billiard paths then we set
$$
\length(p):= (\Length(p_1),..., \Length(p_m)).
$$

\begin{defn}
\label{LS}
A PL path $p: [0,1]\to A$ in $A$ is said to be an LS path with respect to the root
system $R$ if:

1. $p(0), p(1)\in P(R)$.

2. For each $t\in [0,1]$ we have
$$
p'_-(t)\ge_{W_{p(t)}} p'_+(t)
$$

3. There is a $W_{p(t)}$--chain from $p'_-(t)$ to $ p'_+(t)$ which is maximal as a
$W$-chain when we regard $p'_-(t)$ and $ p'_+(t)$
as vectors in $V$.
\end{defn}

\begin{rem}
This definition is a slight generalization of the Littelmann's definition in \cite{Littelmann2},
where it is  assumed that $p(0)=0$.
\end{rem}

\begin{defn}
A PL path $p: [0,1]\to A$ is said to be a {\em Hecke path} if it satisfies properties 1 and 2 
in Definition \ref{LS}. 
\end{defn}

Note that each Hecke path is necessarily a billiard path. 
Let $t_i, i=1,...,n,$ denote the {\em break-points} of $p$, i.e. such that $p$ is
not geodesic at $t$. We then obtain a chain
$$
p'(0) \ge p'_+(t_1)\ge ...\ge p'_+(t_{n-1}) \ge p'(1).
$$
Using property (3) we extend this chain to a maximal $W$-chain
$$
\eta_0=p'(0)\ge ... \ge \eta_m =p'(1).
$$
Then the translated path $p-p(0)$ is the concatenation of the geodesic paths
$$
\pi_{a_i \eta_i}, a_i\ge 0, i=1,...,n.
$$

We also need a generalization of the concept of an LS path described below:

\begin{defn}
[\cite{KM}]
Suppose that $p_1,...,p_m: [0,1]\to A$ are LS paths. Their concatenation
$p_1* ...* p_m$ is called a generalized LS path if for each $i=1,...,m-1$ we have:
$$
(p_i')_{-}(1)\gtrsim (p_i')_{+}(0).
$$
\end{defn}

We will use the notation $LS(R)$ and $LS_1(R)$ to denote the sets of LS paths and generalized
LS paths with respect to the root system $R$.
In fact, in this paper we will be using only (generalized) LS paths
$p=p_1* ...* p_m$ such that
$$
\length(p)= (n_1 \varpi_1,..., n_m\varpi_m), n_i\in \Z_+,
$$
where $\varpi_i$ are fundamental weights of $R$.

\subsection{The path model}

The following theorem is a version of Littelmann's rule for decomposing tensor products:

\begin{thm}
[\cite{KM}, Corollary 5.22] \label{Lit2} The tensor product $V_\la\otimes V_\mu$
contains $V_\nu$ as a subrepresentation if and only if there
exists a generalized LS path $p$ with $p(0)=\la, p(1)=\nu$ so that

1. $$
\length(p)=(\mu_1,...,\mu_m),
$$
2. $$
\mu_i\in \Z_+\varpi_i, \forall i=1,...,m,
$$
3. $$
\sum_{i=1}^m \mu_i=\mu
$$
4. $p$ is contained in $\De$.
\end{thm}

%\begin{thm}
%(P. Littelmann, \cite{Littelmann2}) $\mu\in Weight(V_\la)$ if and
%only if there exists an LS path $p$ with $\Length(p)=\la$ such
%that $p(0)=0, p(1)=\mu$.
%\end{thm}

The following lemma easily follows from the above theorem, one can also derive it directly
from the definition of a generalized LS path:

\begin{lem}
\label{trivial}
If $p$ is a generalized LS path then
$$p(0)-p(1)- \Length(p)\in Q(R).$$
\end{lem}

\subsection{The saturation theorem}

In this section we discuss the {\em
Saturation Theorem} proven in \cite{KM}.

Recall that in \cite{KLM3} we have defined two constants $k_R$ and $k_w$ associated with 
the root system $R$. The {\em saturation constant} $k_R$ is defined by the property that 
it is the least integer $k\ge 1$ such that for each vertex $v$ of the Euclidean 
Coxeter complex $(A, W_{R,aff})$ we have:
$$
kv\in P(R),
$$
i.e. is a special vertex. The constant $k_w$ is the least integer
$k\ge 1$ with the following property:

Suppose that $F$ is a face of $(A, W_{R,aff})$ invariant under an
isometry of $A$ and let $v$ denote the barycenter of $F$. Then
$kv$ is a special vertex.

We have proven in \cite{KLM3} that for the root systems $B_\ell$
and $C_\ell$ we have:
$$
k_R=k_w=2.
$$
In particular, $k_R$ and $k_w$ do not change if we replace $R$ with $R^\vee$.

Let $\ul{G}$ be a reductive algebraic group over $\Z$; set
$G=\ul{G}(\C)$. For   a nonarchimedian local field $\K$
(e.g. $\K=\Q_p$) we let $X$ denote the Bruhat-Tits building associated with
$\ul{G}^\vee(\K)$, where $\ul{G}^\vee$ is the Langlands' dual of $\ul{G}$.

\begin{thm}
\label{saturation}
1. Suppose that $\si=(\la,\mu,\nu)\in P(R)^3\cap \P(G)$ is such
that $\la+\mu+\nu\in Q(R)$. Then
$$
k_R^2\cdot \si \in Tens(G).
$$
2. Suppose that $\si=(\la,\mu,\nu)\in P(R)^3\cap \P(G)$. Then
$$
k_R k_w\cdot \si \in Tens(G).
$$

\end{thm}
 \proof The first assertion is the {\em
Saturation Theorem} 1.8 of \cite{KM}. We prove the second assertion.

Then, since $\si \in \P(G)=D_3(X)$, there
exists a geodesic triangle $[z, x, y]$ in $X$ with the $\De$-side
lengths $(\la,\mu,\nu)$ (see \cite{KLM2}). According to Theorem 7.16
of \cite{KLM3}, there exists a geodesic triangle $[z', x',
y']\subset X$ whose verticies are special verticies of $X$ and whose
$\De$-side lengths are $(\la', \mu', \nu'):=(k_w\la, k_w\mu,
k_w\nu)$. Therefore, by Part 2 of Theorem 1.8 in \cite{KM},
$$
k_R(\la', \mu', \nu')\in Tens(G). \qed
$$

In section \ref{B2C2} we will need the following improvement of Theorem
\ref{saturation} in the case of the root system $R=B_2\cong C_2$.

\begin{thm}
\label{B2saturation}
Suppose that $G$ has the root system $R\cong B_2$,  
$\si=(\la,\mu,\nu)\in P(R)^3\cap \P(G)$ and $\la+\mu+\nu\in Q(R)$. Then:

1. 
$$
2\cdot \si \in Tens(G).
$$
2. Moreover, there exists a generalized LS path $p\in LS_1(R)$
contained in $\De$, connecting $2\la$ to $2\nu=2\nu^*$, so that $\length(p)=(2\mu_1,2\mu_2)$, 
$\mu_1+\mu_2=\mu$, $\mu_i\in \Z_+\varpi_i$ and all break-points of $p$, 
with possible exception of ones occurring on the boundary of
$\De$, are special verticies.
\end{thm}
\proof The first assertion follows from the second. 
%is  Proposition 8.34 of \cite{KLM3}. 
The proof of the second assertion is a
variation on the proof of the Saturation Theorem 1.8 given in \cite{KM} 
so here we will give only a sketch and refer the reader to \cite{KM} for the details.

Let $X$ be the Euclidean building as above. 
Then the assumptions that $\si=(\la,\mu,\nu)\in P(R)^3\cap \P(G)$ and
$\la+\mu+\nu\in Q(R)$ imply that
there exists a geodesic triangle $[\t{x}, \t{y}, \t{z}]\subset X$ whose verticies are verticies of $X$
and whose $\De$--side-lengths are $\mu, \nu, \la$, see \cite[Corollary 7.12, Part 1]{KLM3}. 
We let $\t{A}\subset X$ denote an apartment containing
the segment $\ol{\t{x} \t{y}}$ and let $\t{\De}\subset \t{A}$ denote a translate of a Weyl chamber in $\t{A}$,
so that the tip of $\t\De$ is at $\t{x}$, and $\ol{\t{x} \t{y}}\subset \t{\De}$.
We identify the fundamental weights $\varpi_1, \varpi_2$ with vectors in $\t{A}$ so that
$\mu=m_1\mu_1+ m_2\mu_2$, where $\mu_i$ are multiples of $\varpi_i$, $i=1, 2$.
We then replace the segment
$\ol{\t{x} \t{y}}$ with the concatenation
$$
\t{p}= \pi_{\mu_1}* \pi_{\mu_2}\subset \t{\De}
$$
The path $\t{p}$ connects $\t{x}$ to $\t{y}$. 
%We leave to the reader to verify the following elementary
%observation:

\begin{lemma}
\label{obs}
Let $h: A\to A$ be a dilation by $2$ which fixes a special vertex. Then 
the path $\t{p}$ crosses walls of $\t{A}$ transversally only at points $v$ 
such that $h(v)$ are special verticies.
\end{lemma}
\proof To simplify the notation we identify the apartment $\t{A}$ with the model apartment $(A, W_{aff})$. 
Suppose first that $\t{x}$ is a special vertex.  Since $\mu_1\in P(R)$, the end-point  
of the path $\pi_{\mu_1}$ is also a special vertex. Thus the paths $\pi_{\mu_1}, \pi_{\mu_2}$ (and therefore $\t{p}$) 
are entirely contained in the 1-dimensional skeleton of the simplicial complex $(A, W_{aff})$. Therefore these paths  
cross walls transversally only at the verticies of this complex. However $k_R=2$ for $R=C_2$ means that 
for each vertex $v\in (A, W_{aff})$, its image under dilation $h(v)$ is a special vertex. 
Hence the claim follows in this case. 

\begin{figure}[tbh]
\centerline{\epsfxsize=4.5in \epsfbox{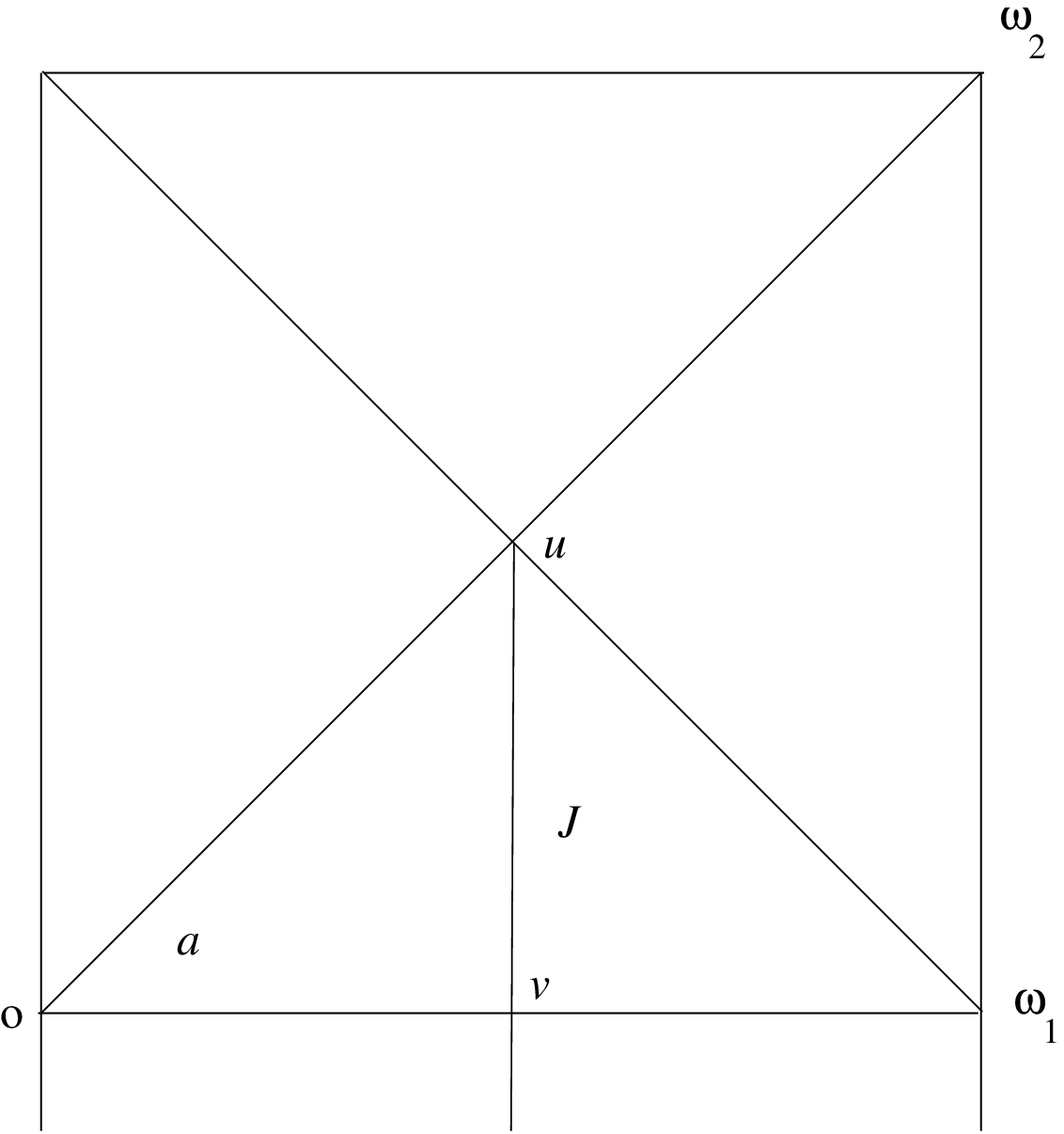}}
\caption{\sl }
\label{special.fig}
\end{figure}

Suppose that $\t{x}$ is not a special vertex. Nevertheless, this point is a vertex of an alcove $a\subset 
(A, W_{aff})$. The break-point $\t{u}$ of the path $\t{p}$ is also a vertex of $(A, W_{aff})$ 
which has the same type as $\t{x}$. Consider now a pair of points $u\in a, s\in A$ which are nonspecial verticies  
so that the geodesic segment $J=\ol{us}\subset (A, W_{aff})$ is parallel to an element of $W\cdot \varpi_1$ or 
$W\cdot \varpi_2$. Then for each point $v$ of transversal intersection of $J$ with walls of $(A, W_{aff})$ we have:

\begin{itemize}
\item Either the type $\theta(v)$ of $v$ is a vertex of $a$, in the case when $J$ is parallel to an element of 
$W\cdot \varpi_2$,

\item Or the type $\theta(v)$ equals $\frac{1}{2}\varpi_2$, in the case when $J$ is parallel to an 
element of $W\cdot \varpi_1$, see Figure \ref{special.fig}.
\end{itemize}

\noindent In either case, $h(v)$ is again a special vertex of $(A, W_{aff})$. \qed

\begin{rem}
The direct generalization of Lemma \ref{obs} fails in the case of the root system $B_3$ where 
the path $\t{p}$ is the concatenation of three geodesic segments parallel to the fundamental weights. 
\end{rem}

Recall that in \cite{KM} we have defined a family of projections
$f=Fold_{z,k,\De}: X\to \De$, where $\De\subset A$ is a Weyl chamber with
tip at $o$. Here $k\in \N$, $z$ is a vertex in $A$. The mapping
$f$ is the composition of three maps:

1. A retraction $Fold_{a,A}: X\to A$ with respect to an alcove $a$
containing $z$.

2. Dilation $h: A\to A$ by $k$ so that $h(z)=o$.

3. Projection ${\mathbb P}_\De: A\to \De$.

Then
$$
f= {\mathbb P}_\De \circ h \circ Fold_{a,A}.
$$

Consider the path
$$
p:= Fold_{z, 2,\De}(\t{p}),
$$
where $z:=\t{z}$. It was shown in \cite{KM} that the
path $p$ is in $LS_1(R)$ and
$$
\length(p)= 2\length(\t{p})=2(\mu_1,\mu_2).
$$
It is clear that $p$ is contained in $\De$ and that this path connects $x:=f( \t{x})$ to $y:=f(\t{y})$
where $\ov{ox}=2\la, \ov{oy}=2\nu$. Observation \ref{obs} implies that all break-points of the
path $q:=h\circ Fold_{a,A}(\t{p})$ are special verticies of $A$. Their images  under the projection
${\mathbb P}_\De$ are also special. The projection ${\mathbb P}_\De$ may introduce new break-points
in the path $p$ (i.e. break-points which are not projections of break-points of $q$). However
such points necessarily belong to the walls of $\De$. \qed

\begin{rem}
It is proven in Proposition 8.34 of \cite{KLM3} that in Part 1 of Theorem \ref{B2saturation}, 
the assumption that $\la+\mu+\nu\in Q(R)$ can be omitted. Theorem \ref{mainBC} 
proven in section \ref{sect:BC} of the present paper, 
provides an alternative proof of this result.  
\end{rem}

\section{Decomposition of $Tens(G)$ into elementary subsets}

\begin{defn}
Call a subset $E\subset \Z^n$ {\em elementary} if it is defined
via a finite system of (non-strict, possibly inhomogeneous) linear
inequalities with rational coefficients and congruences, i.e.
equations of the form $h(x)\in \Z$, where $h$ is a rational
linear function.
\end{defn}

Note that each system of congruence conditions on $x$ is equivalent to the
requirement that $x$ belongs to a coset of a subgroup in $\Z^n$. By adding a linear
equation to the system of inequalities, we  can reduce a system of congruence conditions on $x$ to
the requirement that $x\in L+z$, where $L$ is a sublattice in $\Z^n$
(i.e. a rank $n$ subgroup) and $z$ is a certain element of $\Z^n$.

The next proposition follows for instance from \cite[Theorem 1]{C}
(we are grateful to C.~Laskowski for this reference), but we give
an elementary proof for the sake of completeness.

\begin{prop}\label{logic}
If $E\subset \Z^n$ is an elementary subset then its projection to
$\Z^{n-1}$ is a finite union of elementary subsets.
\end{prop}
\proof We let $(x,a)$ denote coordinates in $\Z^n$ so that $x$ is
the coordinate in $\Z^{n-1}$. Suppose that the elementary set $E$
is given by the linear inequalities
$$
F(x,a) \in (\R_+)^m
$$
and congruences $(x,a)\in L$, where $L$ is a translate of a
sublattice in $\Z^n$. Therefore, up to changing coordinates via an
integer translation, $L$ is a sublattice in $\Z^n$. Since $L$ has
finite index in $\Z^n$, there exists an integer $\kappa$ so that
$\kappa\cdot \Z^n \subset L$; therefore $L$ is a union of finitely
many cosets
$$
L_i:= z_i+ \kappa\cdot \Z^n.
$$
By restricting to the cosets $L_i$ and making the linear changes
of variables $(x,a)\to (x,a)- z_i$, we reduce the proof to the
case when $L$ has the form $\kappa\cdot \Z^n$, which we assume
from now on.

Let $\t{P}$ denote the convex polyhedron $\{(x,a): F(x,a) \in
(\R_+)^m\}$. Let $P$ denote the projection of $\t{P}$ to
$\R^{m-1}$ under the map $p: (x,a)\mapsto x$; this set is again a
convex polyhedron. Then we can subdivide $P$ into a finite union
of convex polyhedra $P_i$ (each given by a linear system of
inequalities with rational coefficients) such that for each $i$
the set $p^{-1}(P_i)\cap \t{P}$ is given by the two inequalities
$$
g_i(x)\le a\le f_i(x)
$$
where $g_i, f_i$ are linear functions with rational coefficients.
From now on we fix $i$ and set $f:= f_i, g:= g_i$: We will show
that $p(E)\cap P_i$ is a finite union of elementary sets.

Observe that the restriction of $f/\kappa$ to $\kappa\cdot
\Z^{n-1}$ takes only finitely many values (mod
$\Z$)
$$
t_j\in [0, 1), j=1,...,J.$$

Let $\La_j$ denote the coset in $\kappa\Z^{n-1}$ such that the
value of $f/\kappa$ on $\La_j$ equals $t_j$ mod $\Z$. Let
$[f(x)/\kappa]$ denote the integer part. Then the condition that
there exists $a\in \kappa\Z$ such that
$$
g(x)\le a\le f(x)
$$
is equivalent to
$$
g(x)/\kappa\le [f(x)/\kappa]
$$
i.e.
$$
g(x)/\kappa\le f(x)/\kappa - t_j, x\in \La_j.
$$

Therefore the projection of $E\cap p^{-1}(P_i)$ to $P_i$ equals
$$
E_{ij}=\bigcup_{j=1}^{J} \{ x\in P_i\cap \La_j : g(x)\le f(x)-
\kappa t_j\}.
$$
It is clear that each $E_{ij}$ is an elementary set. \qed

\medskip
Our next goal is to show that the semigroup $Tens:= Tens(G)$ is a finite union of elementary
sets, where $G$ is a complex reductive Lie group. Let $L=X^*(T)$ denote the cocharacter lattice
of a maximal torus $T\subset G$. We let $\varpi_1,...,\varpi_\ell$ denote the fundamental weights
of $G$, $\Pi:= W\cdot \{\varpi_1,...,\varpi_\ell\}$, where
$W$ is the Weyl group of $G$. Let $\al_1,...,\al_\ell$ denote the simple roots in $R$.
Let $\al_i^\vee$ denote the coroots.

For each $\eta\in \Pi$ let $\iota(\eta)$ denote the number $i$ such that $\eta\in W\cdot
\varpi_i$. Let $\Si$ denote the set of {\em generalized chains}
$\si$ in $\Pi$, i.e. sequences
$$
\eta_0 \gtrsim \eta_1 \gtrsim ... \gtrsim \eta_m
$$
of elements of $\Pi$ so that
$$
\iota(\eta_0)\le \iota(\eta_1)\le ... \le \iota(\eta_m).
$$
Let $\Si_{max}$ denote the collection of {\em maximal} generalized
chains as above, i.e. chains of maximal length from $\eta_0$ to
$\eta_m$.

Let $\la, \mu, \nu\in \De\cap L$ be dominant characters,
$$
\mu= \sum_{i=1}^\ell n_i \varpi_i.
$$
We consider broken geodesic paths $\pi: [0,1]\to V=P(R)\otimes \R$
modeled on the chain $\si$ as above, with $p(0)=0$, i.e. concatenations of paths
$$
\pi_{a_i \eta_i}, i=0,...,m,
$$
where $a_i\in \frac{1}{k_R^2} \Z$, $a_i\ge 0$. We require that
$$
\length(\pi)= (n_1\varpi_1,...,n_\ell\varpi_\ell).
$$
This means that
$$
a_1+...+a_{J_1}= n_1, a_{J_1+1}+...+a_{J_2+2}= n_2,...
$$

We define the partial sums
$$
S_r=\sum_{j=1}^r a_j \eta_j
$$
for $1\le r\le m$. Observe that $\iota(\eta_r)=\iota(\eta_{r+1})$ if nd only if there
exists a reflection $\tau_r=\tau_{\al_r^\vee}$ which carries $\eta_r$ to $\eta_{r+1}$. 
Then the path $\pi$ is a generalized LS path if and only if the
following condition is satisfied:

For each $r$ such that $\iota(\eta_r)=\iota(\eta_{r+1})$
we have
$$
\al_r^\vee(S_r)= \sum_{j=1}^r a_j \al_r^\vee(\eta_j) \in \Z.
$$
Note that these conditions imply that if $\iota(\eta_r)\ne
\iota(\eta_{r+1})$ then $S_r$ is necessarily a special vertex,
i.e. all coroots take integer values at this point.

Since $\al_r^\vee(\eta_j) \in \Z$ for all $r, j$, it follows that
the above integrality condition is a {\em congruence} condition on
the $m$-tuple $a=(a_0,...,a_m)$.

According to Theorem \ref{Lit2}, we have: 
$(\la,\mu,\nu)\in Tens(G)$ if and only if there exists a generalized
LS path $\pi$ as above so that:

(1) $\la+\pi(1)=\nu^*$.

(2) The entire path $\la+\pi(t)$ is contained in the positive
chamber $\De$, i.e. for each simple coroot $\al^\vee$ and each
partial sum
$$
S_r=\sum_{j=1}^r a_j \eta_j
$$
we have:
$$
\al^\vee(\la+S_r)\ge 0.
$$

Let $\la= \sum_{i=1}^\ell m_i \varpi_i$, $\nu= \sum_{i=1}^\ell k_i
\varpi_i$. Set
$$
x:= (m_1,...,m_\ell, n_1,...,n_\ell, k_1,..., k_\ell).
$$

Note that the condition (1) is a linear equation with integer
coefficients on the vector $a$ and the condition (2) has the form
of a system of linear inequalities with integer coefficients.
Therefore for each  generalized chain $\si$ as above, the set
$$
E_\si:=\{(x, a) : \la+p(1)=\ga^*, \al^\vee(\la+S_r)\ge 0, \forall \al\in \Phi, \hbox{and~}
 \al_r^\vee(S_r) \in \Z, r=0,...,m-1 \}
$$
is an elementary set.

\begin{rem}
Instead of using path model in the above argument one can use the polytopal model from \cite{BZ}.
\end{rem}

Consider the projection $p(E_\si)$ of $E_\si$ to the
$x$-coordinate. By applying Proposition \ref{logic} inductively we
conclude that $p(E_\si)$ is a finite union of elementary sets.
Therefore the union
$$
Tens(G)= \bigcup_{\si\in \Si_{max}} p(E_\si)
$$
is also a finite union of elementary sets.

Thus we have proved the following analogue of Laskowski's
theorem in \cite{Laskowski}:

\begin{thm}
The semigroup $Tens(G)$ is a finite union of elementary sets.
\end{thm}

\section{Deep subcone}
\label{subcone}

The goal of this section is to show that deep inside of the cone
$D_3=\P(G)$ there is a subcone of the form $D_3+\si$ (for some
$\si\in D_3$), such that
$$
(D_3+\si) \cap \La\subset Tens(G).$$

We will also present an explicit computation of this subcone for
the root system $B_\ell$. In what follows we normalize the roots so that the Euclidean norm of each
coroot is either $1$ or $2$.

Given an irreducible representation $V_\la$ let $Weight(V_\la)\subset P(R)$ denote the set of weights of
 $V_\la$. Define the set
$$
S_3(G)^{1,2}:= \{ (\la,\mu,\nu)\in \De^3: \la \gg \mu, \nu^*= \la+\be, \hbox{for some~} \be\in Weight(V_\mu)\}
$$
Here $\la\gg \mu$ iff $\la + \be \in \De$ for all $\be\in
Weight(V_\mu)$.

Below we give a more explicit description of the
above subset $S_3(G)^{1,2}$ in terms of linear inequalities and
congruences.

We will be using the following notation: If $C\subset V$ is a convex subset then $\la\ge_C \mu$ iff $\mu-\la\in C$.
We use the notation $\De^*$ for the convex cone
{\em dual to $\De$}, i.e.
$$
\De^*=\{v\in V: v\cdot u\ge 0, \forall u\in V\}.
$$

Define the lattice
$$\La=\{(\mu, \nu, \la): \mu, \nu, \la\in  L, \la+\mu+\nu\in Q(R)\},$$
where $L=X^*(T)$. We let $C^{1,2}$ be convex polyhedral cone in $\De^3$ given by the
following inequalities:

1.
$$
 w\mu^* \le_{\De} \la,\quad  \forall w\in W.$$

2.
$$
w\nu^* \le_{\De^*} w\la +\mu ,\quad \forall w\in W.
$$

\begin{prop}
\label{prel}
$$S_3(G)^{1,2}= C^{1,2} \cap \La. $$
\end{prop}
\proof

\begin{lem}\label{La}
$$
\la\gg \mu \iff  w\mu^* \le_{\De} \la,\quad  \forall w\in W.
$$
\end{lem}
\proof
First let's check that if $\la+w\mu\in \De$ for all $w\in W$ then
$$
\la+\be \in \De, \forall \be\in Weight(V_\mu).
$$
Indeed, let $\be\in Weight(V_\mu)$. We write $\be$ as a convex combination of extreme weights $\{w\mu, w\in W\}$:
$$
\be= \sum_{w\in W} t_w (w\mu), \sum_{w\in W} t_w=1, 0\le t_w\le 1.
$$
Then
$$
\la+\be= \la+ \sum_{w\in W} t_w (w\mu)= \sum_{w\in W} t_w \la+ \sum_{w\in W} t_w (w\mu)=\sum_{w\in W} t_w (\la+ w\mu).
$$
By assumption, $\la+w\mu\in \De$. Since $\De$ is convex, $\la+\be\in \De$ as well. Therefore
$$
\la\gg \mu \iff \la+ w\mu \ge_\De 0, \forall w\in W \iff  \la+ ww_0\mu \ge_\De 0, \forall w\in W,
$$
where $w_0\in W$ is the longest element.
$$
\la+ ww_0\mu \ge_\De 0  \iff ww_0 (-\mu) \le_\De \la \iff w\mu^* \le_\De \la,
$$
since $\mu^*= w_0 (-\mu)$. \qed

We refer the reader to \cite{H} for the proof of the following:

\begin{lem}
\label{Lm}
For a dominant weight $\mu$ and a weight $\be$ we have
$$
\be \in Weight(V_\mu)\iff w\be \le_{\De^*} \mu, ~~ \forall w\in W, \hbox{~~and~~} \mu-\be \in Q(R).
$$
\end{lem}

\begin{lem}
\label{Lb}
Let $\la, \mu, \nu$ be dominant weights. Then
$$
\nu^*=\la+\be, \hbox{~for some~} \be\in Weight(V_\mu)
$$
if and only if
$$
w\nu^*\le_{\De^*} w\la +\mu, \hbox{~~and~~} \mu-\nu^*+\la\in Q(R).
$$
\end{lem}
\proof
$$
\nu^* - \la\in Weight(V_\mu) \iff w(\nu^* - \la) \le_{\De^*} \mu, \forall w\in W,  \hbox{~~and~~} \mu-(\nu^*-\la)\in Q(R).
$$
\qed

It remains to prove

\begin{lem}
\label{Lc}
$$\mu-\nu^* +\la\in Q(R) \iff \mu+\nu+\la\in Q(R).$$
\end{lem}
\proof It suffices to prove that $\mu-\nu^* +\la-(\mu+\nu+\la)\in Q(R)$, equivalently,
$$
\nu+\nu^*= \nu - w_0\nu \in Q(R).
$$
However for each $\nu \in P(R), w\in W$ we have $w\nu\in Weight(V_\nu)$, therefore
$$
\nu - w\nu\in Q(R),
$$
by Lemma \ref{Lm}. \qed

This concludes the proof of Proposition \ref{prel}. \qed

The following result is standard, we are grateful to T.~Haines and S.~Kumar for pointing out this result to us and
explaining the proofs:

\begin{thm}
\label{kumar}
 The subset $S_3(X)^{1,2}$ is contained in $Tens(G)$.
\end{thm}
\proof We will present the proof of this result using Littelmann's
path model.

We need the following

\begin{lem}
Suppose that $p$ is an LS path of $\De$-length $\mu$. Then the
path $p$ is entirely contained in the convex hull of the $W$-orbit
$S:=W(\mu)$.
\end{lem}
\proof By definition of an LS path, for each $t\in [0,1]$ there exist $t_i\ge 0$ and $\nu_i\in W\mu$ such that
$$
t=t_1+...+t_j
$$
and
$$
p(t)= t_1\nu_1+...+t_j\nu_j.
$$
Therefore, since $t\in [0,1]$, and the convex hull of $S$ contains the origin,
the subconvex combination $p(t)$ is contained in the convex hull of $S$. \qed

\medskip
Suppose now that $\si=(\la,\mu,\nu) \in S_3(X)^{1,2}$. Then
$\nu^*=\la+\be$ for some weight vector $\be$ of the representation
$V_\mu$. Then, according to \cite{Littelmann2}, there exists an LS
path $p$ of the $\De$-length $\mu$ such that $p(0)=0, p(1)=\be$.
By the above Lemma, the path $p$ is entirely contained in the
convex hull of $W \mu$. Consider the path $q(t):= \la+ p(t)$.
We claim that this path is entirely contained in $\De$. Indeed,
$$
q(t)\in \la+Hull(W \mu).
$$
Since $\la\gg \mu$, for each vector $\ga\in W\cdot \mu$, we have:
$$
\la+\ga\in \De.
$$
Thus $\la+Hull(W \mu)\subset \De$. On the other hand,
$q(1)=\la+\be=\nu^*$. Therefore, according to Theorem \ref{Lit2}, $(\la,\mu,\nu)\in Tens(G)$. \qed

 We next observe that the cone $C^{1,2}$ has
nonempty interior. Indeed, first choose $\mu\in Int(\De)$. Then
take $\la\in \De$ such that
$$
d(\la, \D \De) > |\mu|.
$$
Finally, pick $\nu^*$ sufficiently close to $\la$ so that
$$
|\la-\nu^*|< d(\mu, \D \De^*).
$$

Any triple $(\la,\mu,\nu)$ chosen like this satisfies the strict inequalities
\begin{equation}
 w\mu^* <_{\De} \la ,\quad  \forall w\in W,\end{equation}

\begin{equation}
w\nu^* <_{\De^*} w\la +\mu ,\quad  \forall w\in W.
\end{equation}
and therefore belongs to the interior of $C^{1,2}$. Our next goal is to apply the above observations to show that
the cone $\P(G)$ contains a subcone $\si+\P(G)$ such that
$$
\si+\P(G)\cap \La= Tens(G)\cap \La.
$$

Since $C^{1,2}$ is a homogeneous cone with nonempty interior, it
contains metric balls $B(\si_0, R)$ of arbitrarily large radius
$R$. Let $k:= k_R k_w$. Choose $R$ larger than the diameter of a
fundamental domain $F$ for the lattice $k\cdot P(R)$. (Here $F$ is
a certain fundamental parallelepiped containing the origin.)
Without loss of generality we assume that the ball $B(\si_0, R)$ is centered at
a point $\si_0\in k\cdot P(R)$.

Suppose now that $\tau\in \La \cap D_3(X)$. Then there exists a
point $\si\in k\cdot P(R)$ such that $\tau\in F+\si$. Let
$\kappa:= \si-\si_0$. Since
$$
F+\si= (F+\si_0) + \kappa,
$$
there exists a point $\tau_0\in F+\si_0$ such that
$\tau=\tau_0+\kappa$.

Note that $k\cdot P(R)\subset Q(R)$, since the index $|P(R):
Q(R)|$ divides $k_R k_w$, see \cite{KLM3}, Table 11. Therefore
$$
\si_0\in k\cdot P(R)^3\subset \La
$$
and hence $\tau_0=\tau-\kappa= \tau-\si+\si_0 \in \La$. Since
$diam(F)\le R$, $F+\si_0$ is contained in $B(\si_0, R)\subset
C^{1,2}$. Therefore, by  Theorem \ref{kumar}, $\tau_0\in  F+\si_0
\subset C^{1,2} \subset Tens(G)$. By the triangle inequality, if
$$
d(\tau, \D D_3(X))\ge |\si_0| +diam(F),
$$
then $\kappa=\si-\si_0$ belongs to $\P(G)$. On the other hand,
since $\si, \si_0\in k\cdot P(R)$, it follows that $\kappa\in
k\cdot P(R)$ and therefore, by the saturation theorem,
$$
\kappa\in Tens( G).
$$
Therefore, since $\tau_0\in Tens(G)$ and $Tens(G)$ is a semigroup, it follows that
$$
\tau=\tau_0+\kappa
$$
also belongs to  $Tens(G)$. Hence we have proven

\begin{thm}
\label{deep} Suppose that $\tau\in \P(G)\cap \La$ is such that
$d(\tau, \partial \P(G))\ge |\si_0| +diam(F)$. Then $\tau\in
Tens(G)$.
\end{thm}

Next, the linear inequalities defining $\P(G)$ have the form
$$
\psi_j(x)  \ge 0,
$$
where either $|\psi_j|\le 2$ (in the case of the inequalities
$\psi_j\ge 0$ defining the chamber $\De$) or
$|\psi_j|=\sqrt{3}|\varpi_{i_j}|$ for a certain fundamental weight
$\varpi_{i_j}$ (in the case of the stability inequalities). Here
we are using the Euclidean norm of linear functionals. Set
$m:=\max_i |\varpi_i|$ and note that, because of our normalization of the lengths of the coroots,
$m \ge 2$. Therefore,
$$
\max_j |\psi_j|=m.
$$
Suppose that $\tau-\kappa\in \D \P(G)$. Then there exists some $\psi_j$ such that $\psi_j(\tau-\kappa)=0$, i.e.
$$
|\psi_j(\tau)|= |\psi_j(\kappa)|\le |\psi_j|\cdot |\kappa| \le m|\kappa|
$$
and $|\kappa|\ge |\psi_j(\tau)|/m$. Thus, if
$$
\psi_j(\tau)  \ge m(|\si_0| +diam(F)), \forall j,
$$
then $d(\tau, \partial \P(G))\ge |\si_0| +diam(F)$. Therefore, if we define an inhomogeneous subcone $\P'(G)$
in $\P(G)$ by the linear inequalities
$$
\psi_j(\tau)  \ge m(|\si_0| +diam(F)), \forall j,
$$
we obtain:

\begin{cor}
\label{mainc}
The intersection $\P'(G)\cap \La$ is contained in $Tens(G)$.
\end{cor}

Below we make an explicit computation for the root system
$R=B_\ell$ (i.e. when $G=Sp(2\ell, \C)$). We use the Bourbaki
coordinates \cite{Bourbaki},
to describe this root system. Since for this root system $w_0=-1$, it follows that $\la=\la^*$
for each $\la\in \De$. We let $\{\eps_i\}$ denote the
standard orthonormal basis in $V$. Then for each $i$,
$\varpi_i=\sum_{j=1}^i \eps_i$. Moreover, the simple roots are
$\al_i=\eps_{i-1}-\eps_i$, $i<\ell$ and
$\al_\ell=\eps_\ell$. The positive chamber $\De$ is
given by the inequalities $\al_i\cdot v\ge 0, i=1,...,\ell$.

\begin{thm}
Suppose that $R=B_\ell$ ($\ell\ge 2$), and $\tau\in \La$ is such
that for each linear functional $\psi$ which appears in the system
of stability inequalities and chamber inequalities we have:
\begin{equation}
\label{0} \psi (\tau)  \ge 2\ell^2 (\ell +1)(4\ell+5) +6\ell.
\end{equation}
Then $\tau\in Tens(Sp(2\ell, \C))$.
\end{thm}
\proof We will use the notation from the proof of Theorem \ref{deep}.
In order for $\tau$ to be in  $Tens(Sp(2\ell, \C))$ we need two things:

\begin{equation}
\label{1}
\si_0+F\subset C^{1,2}
\end{equation}
and
\begin{equation}
\label{2}
\kappa=\tau-\tau_0\in \P(G), \quad \forall \tau_0\in F+\si_0.
\end{equation}

We  simplify the discussion by observing
that in the case of the root system $B_\ell$ we have $k=k_R=k_w=k_G=2$ (using the notation of \cite{KLM3}),
which means that instead of working with the lattice $\La$ we can work with the lattice
$L^3$, where $L=P(R^\vee)$:

According to Part 2 of Theorem \ref{saturation}, for each $\si\in L^3\cap \P(G)$, the vector $k^2\si$ belongs to $Tens(G)$.

The lattice $L$ is just the integer lattice in $\R^\ell$ (using Bourbaki coordinates).
Therefore we choose the fundamental domain $F$ for the sublattice $(4L)^3$ to be the cube
whose edges have length 4 and which is centered at the origin.

The condition (\ref{2}) would follow from:
$$
\psi_j (\tau ) \ge \psi_j (\tau_0)$$
which in  turn is implied by

\begin{equation}
\label{6}
\psi_j (\tau ) \ge \psi_j (\si_0) +\max_{f\in F} \{ \psi_j (f)\}.
\end{equation}

\begin{rem}
$M=\max_{j, f\in F} \{ \psi_j (f)\}$ equals the maximum of all $\psi_j$'s
on the set of points in $V^3$ with coordinates $\in \{0, \pm 1, \pm 2\}$. Here $V=P(R)\otimes \R$
\end{rem}

In what follows we will use the norm
$$
\|(v_1,...,v_\ell)\|= \sum_{i=1}^\ell |v_i|
$$
for vectors $v\in V$.

For the linear functionals $\psi_j$ which come from the stability
inequalities, the maximum
$$
\max_{f\in F} \psi_j(f)
$$
 does not exceed $6\|\varpi_\ell\|=6\ell$. If the inequality $\psi_j\ge 0$ is one
of the chamber inequalities, then $\psi_j(f)\le 4\le 6\ell$ for
all $f\in F$. Thus $M\le 6\ell$.

Hence for the  root system $B_\ell$ to guarantee (\ref{6}) (and hence (\ref{2})) it suffices to require that
\begin{equation}
\label{2c}
\psi_j (\tau)  \ge \psi_j(\si_0) +6\ell, \forall j.
\end{equation}

To get an explicit estimate we have to choose an appropriate $\si_0=(\la_0, \mu_0, \nu_0)$. Set
$$
\mu_0:= (4\ell,..., 8, 4),
$$
$a:= 4\ell+8, s:= 8\ell+4$ and
$$
\nu_0=\la_0:=
(a+(\ell-1)s,..., a+2s, a+s, a).
$$
Clearly $\la_0, \mu_0, \nu_0 \in \De\cap 4P(R)$.

Observe that
$$
\|\mu_0\|=2\ell(\ell+1),
$$
$$
\|\la_0\|=\|\nu_0\|=a\ell + s \frac{\ell(\ell-1)}{2}= \ell(a+
\frac{s\ell-s}{2})= \ell(4\ell^2+2\ell+6)\le 4\ell(\ell+1)^2.
$$
Therefore for each linear functional $\psi$ of the form
$$
\psi(\la,\mu,\nu)= w_1\la \cdot \varpi_i + w_2\mu \cdot \varpi_i + w_3 \nu \cdot \varpi_i
$$
we obtain:
$$
|\psi(\si_0)|\le |\varpi_i|( \|\la_0\|+ \|\mu_0\|+ \|\nu_0\|)\le \ell(
8\ell(\ell+1)^2 + 2\ell(\ell+1))=2\ell^2(\ell+1)(4\ell+5).
$$
Therefore (\ref{2c}) follows from
\begin{equation}
\label{2d} \psi (\tau)  \ge 2\ell^2(\ell+1)(4\ell+5) +6\ell,
\forall \psi.
\end{equation}

We now consider the condition (\ref{1}), i.e. that each
$$
 (\la,\mu,\nu)= \si_0+f, f\in F,
$$
satisfies
\begin{equation}
\label{1a}
w\mu\le_\De \la, \forall w\in W,
\end{equation}
and
\begin{equation}
\label{1b}
w\nu \le_{\De^*} w\la +\mu, \forall w\in W.
\end{equation}

Let $\phi_{ij}$ denote
the linear functionals of the form
$$
\phi_{ij}(\la, \mu, \nu)= \al_i\cdot (\la- w_j\mu)
$$
where $\al_i$ is a simple root and $w_j\in W$.

Let $\eta_{ij}$ denote the linear functionals of the form
$$
\eta_{ij}(\la, \mu, \nu)=   \varpi_i\cdot ( w_j\la +\mu
-w_j\nu), w_j\in W.
$$

Then to guarantee (\ref{1a}) and (\ref{1b}) we need:
\begin{equation}
\label{3} \phi_{ij}(\si_0)\ge \phi_{ij}(f), \forall f\in F,
\end{equation}
and
\begin{equation}
\label{4} \eta_{ij}(\si_0)\ge \eta_{ij}(f), \forall f\in F
\end{equation}
respectively.

We first consider (\ref{4}) since it is simpler. Since
$\nu_0=\la_0$, the inequality (\ref{4}) reads as
$$
\varpi_i\cdot  \mu_0\ge \eta_{ij}(f), \forall f\in F.
$$
The right hand side of this inequality is $\le 6i$, while the left hand-side equals
$$
2(\ell+i)(\ell+1-i).
$$
It is clear that for $1\le i\le \ell$ and each $\ell \ge 2$ we
have:
$$
2(\ell+i)(\ell+1-i)\ge 2(\ell^2-i^2) + 2(\ell+1)(\ell+i)\ge
6(\ell+i)\ge 6i.
$$
This takes care of the condition (\ref{4}).

Now consider (\ref{3}). First note that $\phi_{ij}(f)\le 8$ for
all $i, j$ and $f\in F$.
Next, for all $w\in W$ we have
$$
\al_i\cdot (w\mu_0)\le 8\ell-4$$
for each $i=1,...,\ell-1$
and
$$
\al_\ell\cdot (w\mu_0)\le 4\ell.
$$

On the other hand,
$$
\phi_{ij}(\si_0)= \al_i\cdot(\la_0-w_j \mu_0)= \al_i\cdot \la_0- \al_i\cdot(w_j \mu_0).
$$
We have:
$$\al_i\cdot \la_0=s, \forall i<\ell$$
and
$$
\al_\ell\cdot \la_0=a.$$
Therefore for $i=1,...,\ell-1$,
$$
\phi_{ij}(\si_0)- \phi_{ij}(f)\ge s-8\ell+4- 8 = (8\ell+4)-8\ell-4=0
$$
and for $i=\ell$
$$
\phi_{ij}(\si_0)- \phi_{ij}(f)\ge a- (4\ell+ 8) =  (4\ell+8)- (4\ell+ 8)=0.
$$

Therefore the condition (\ref{3}) holds. Hence the inequality
(\ref{2d}) implies that $\tau\in Tens(Sp(2\ell, \C))$. \qed

%\medskip
%It appears that the following generalization of Theorem \ref{deep} is true:

\begin{conj}
Suppose that $G$ is a complex semisimple Lie group,
$S$ is a stratum in the boundary of the cone $\P(G)$ which contains a point $(\la,\mu,\nu)$ such that all the vectors
$\la, \mu, \nu$ are regular. Then there exists a subcone $S+\si\subset S$ such that
$$
S+\si \cap \La = (S+\si) \cap  Tens(G)
$$
\end{conj}

%This should follow from Bernstein's observation that there are no nontrivial congruence conditions on
%$Tens\cap \R\cdot (\la,\mu,\nu)$.

\section{Computation of the semigroup $Tens(Sp(4,\C))$}
\label{B2C2}

Let $V$ be $\R^2$ with the coordinates $x, y$ and Euclidean
metric. Let $A$ denote the underlying affine space. Let
$R:=C_2\subset V$ and $R^\vee=B_2\subset V^*=V$ be root systems
with the sets of simple roots equal to $\{ \al_1(x,y)=x-y,
\al_2(x,y)=y\}$ and $\{\al_1^\vee(x,y)=x-y, \al_2^\vee(x,y)=2y\}$
respectively. The fundamental weights of $C_2$ are
$\varpi_1=(1,0), \varpi_2=(1,1)$; the root lattice is
$$
Q(R)=\{(x,y)\in \Z^2: x+y\in 2\Z\}
$$
and the weight lattice is $P(R)=\Z^2$.

Let $W_{aff}:=W_{R^\vee,aff}, W_{aff}^\vee:= W_{R,aff}$ denote 
the affine Weyl groups acting on $A$ corresponding to the root systems 
$R$ and $R^\vee$ respectively. 
Note that $W_{aff}\subset W_{aff}^\vee$ is a subgroup of index 2,
so that all verticies of $(A, W_{aff})$ are special verticies of
$(A, W^\vee_{aff})$. The root lattice $Q(R)$ equals the translation subgroup of $W_{aff}$ 
and the weight lattice $P(R)$ is the normalizer of $W_{aff}$ in the group of Euclidean 
translations of $A$.

Our goal is to compute the semigroup $Tens=Tens(R)$, for $R=C_2$.
We know (from \cite{KLM3}) that each triple $\si=(\la,\mu,\nu)\in
Tens$ satisfies:

1. $\si\in D_3:=D_3(C_2)$.

2. $\si\in P(R)^3$.

3. $\la+\mu+\nu\in Q(R)$.

\noindent It was shown in \cite{KLM3} that these conditions are necessary
but not sufficient for $\si$ to be in $Tens$. The goal of this
chapter to find {\em necessary and sufficient conditions}.

\medskip 
We start by observing that since $\la+\mu+\nu\in Q(R)$, it follows
that at least one of the vectors $\la,\mu,\nu$ belongs to
$Q(R)$. 

\begin{convention}
\label{convent}
Throughout the rest of the paper we will assume that
$\mu\in Q(R)$.
\end{convention}

\begin{notation}
\label{note}
We break the vector $\mu$ as $\mu=\mu_1+\mu_2$, where
$\mu_i=n_i\varpi_i$,
$$
\mu_1= (n_1,0), \mu_2=(n_2, n_2), \mu=(n_1+n_2, n_2).
$$
\end{notation}

Since $\mu, \mu_2\in Q(R)$, it follows that $\mu_1\in Q(R)$ as
well, i.e. $n_1$ is even.

\begin{thm}
\label{mainBC} Suppose that $\si\in D_3(R)\cap P(R)^3$ is such
that $\la+\mu+\nu\in Q(R)$. Then $\si\notin Tens$ if and only if
two of the three vectors $\la, \mu, \nu$ belong to $\Z_+\varpi_2$
and $\la+\mu+\nu\notin 2P(R)$.
\end{thm}
\proof The proof of this theorem occupies the rest of this chapter.
Our strategy is to analyze the geometry of generalized LS paths with respect to
the root system $2R$ and show that all such paths (subject to the condition on
$\la,\mu,\nu$ described in the above theorem) can be transformed to
generalized LS paths with respect to the root system $R$. We first do this in the case 
of {\em singular} paths $p_i\in LS(2R)$ with $\Length(p_i)\in \N\varpi_i$ (section \ref{singular}) 
and then use the results to deal with the general $LS_1(2R)$ paths $p=p_1*p_2$ 
(section \ref{regular}). We are able to carry out this approach largely thanks to the fact that 
the root system $C_2$ is rather small and there are not that many chains formed by elements of 
$W\cdot \varpi_i, i=1, 2$.  This allows us to describe the paths in $LS_1(2R)$ rather explicitely. 
A large supply of  $LS_1(2R)$ paths is given by the Saturation Theorem \ref{B2saturation} which 
is another key tool in the proof.

\medskip
To get an idea how the proof of Theorem \ref{mainBC} might proceed, consider the case when 
$\la, \mu, \nu\in Q(R)$ (i.e. the problem of decomposing tensor products of representations of 
$PSp(4,\C)$). Then $\la, \mu, \nu\in 2P(R^\vee)$ and, according to Proposition 8.34 of \cite{KLM3}, 
$$
(\la, \mu, \nu)\in Tens(R^\vee). 
$$
Therefore there exists a path $p\in LS_1(R^\vee)$ contained in $\De$ connecting $\la$ to $\nu$ and having
$$
\length(p)=(\mu_1, \mu_2). 
$$
If we are lucky, all break-points of $p$ occur at special verticies of $(A, W^\vee_{aff})$ 
(compare Theorem \ref{B2saturation}), hence they are at the verticies of  $(A, W_{aff})$, 
see Corollary \ref{C2}. The key then is to replace $p$ with a new $LS_1(R^\vee)$ path $\hat{p}$ 
contained in $\De$ whose break-points occur at the special verticies of $(A, W_{aff})$, the 
crucial lemmas proving this are Lemmas \ref{K1}, \ref{K2}. 
This would imply that $\hat{p}$ is an $LS_1(R)$ path and $(\la, \mu, \nu)\in Tens(R)$.

\subsection{Analysis of the $LS_1(R)$ paths}

In this section we describe generalized LS paths $p$ with respect
to the root system $R$. Recall that every such path is a
concatenation
$$
p=p_1*p_2, \quad \Length(p_i)\in \N \varpi_i, i=1, 2,
$$
where each $p_i$ is an LS path.
%provided that $\Length(p)\in \N \varpi_i, i=1, 2$.
Let $\tau_i$ denote the reflections in the walls $\{\al_i=0\}$, $i=1, 2$.

Suppose that $p$ is a PL path in $A$ parameterized with the
constant speed.
%an LS path (with respect to $R$)
We say that $t$ (or $p(t)$)  is a {\em break-point} of $p$, if $p$
is not geodesic at $t$.

\begin{defn}
We will call a break-point $p(t)$ {\em a turning point} if the
vectors $p_-'(t)$, $p_+'(t)$ are linearly independent and {\em a
point of backtracking} if $p'_-(t)=-p'_+(t)$. %Observe that at each
%backtracking point the path $p$ satisfies {\em Hecke path
%condition}, while at each turning break point --- it does not.
\end{defn}

Note that for $\eta\in V$, each chain in $W\eta$ can be extended to a chain
which starts in $-\De$ and ends in $\De$. We leave the proof of the following lemma to the reader:

\begin{lem}
Suppose that $\eta\in \{\varpi_1, \varpi_2\}$. Then every maximal
chain in $W(\eta)$ has to be of the form:
$$
\eta_0, \eta_1=\tau_1(\eta_0), \eta_2=\tau_2(\eta_1),....
$$
or
$$
\eta_0, \eta_1=\tau_2(\eta_0), \eta_2=\tau_1(\eta_1),....
$$
\end{lem}

By combining these observations, we see that each chain in
$W(\varpi_1)$ is a subchain of:
$$
\eta_1=-\varpi_1, \eta_2=\tau_2(\eta_1), \eta_3=\tau_1(\eta_2),
\eta_4=\tau_2(\eta_3)=\varpi_1
$$
and each chain in $W(\varpi_2)$ is a subchain of:
$$
\eta_1=-\varpi_2, \eta_2=\tau_1(\eta_1), \eta_3=\tau_2(\eta_2),
\eta_4=\tau_1(\eta_3)=\varpi_2.
$$

\begin{figure}[tbh]
\centerline{\epsfxsize=4.5in \epsfbox{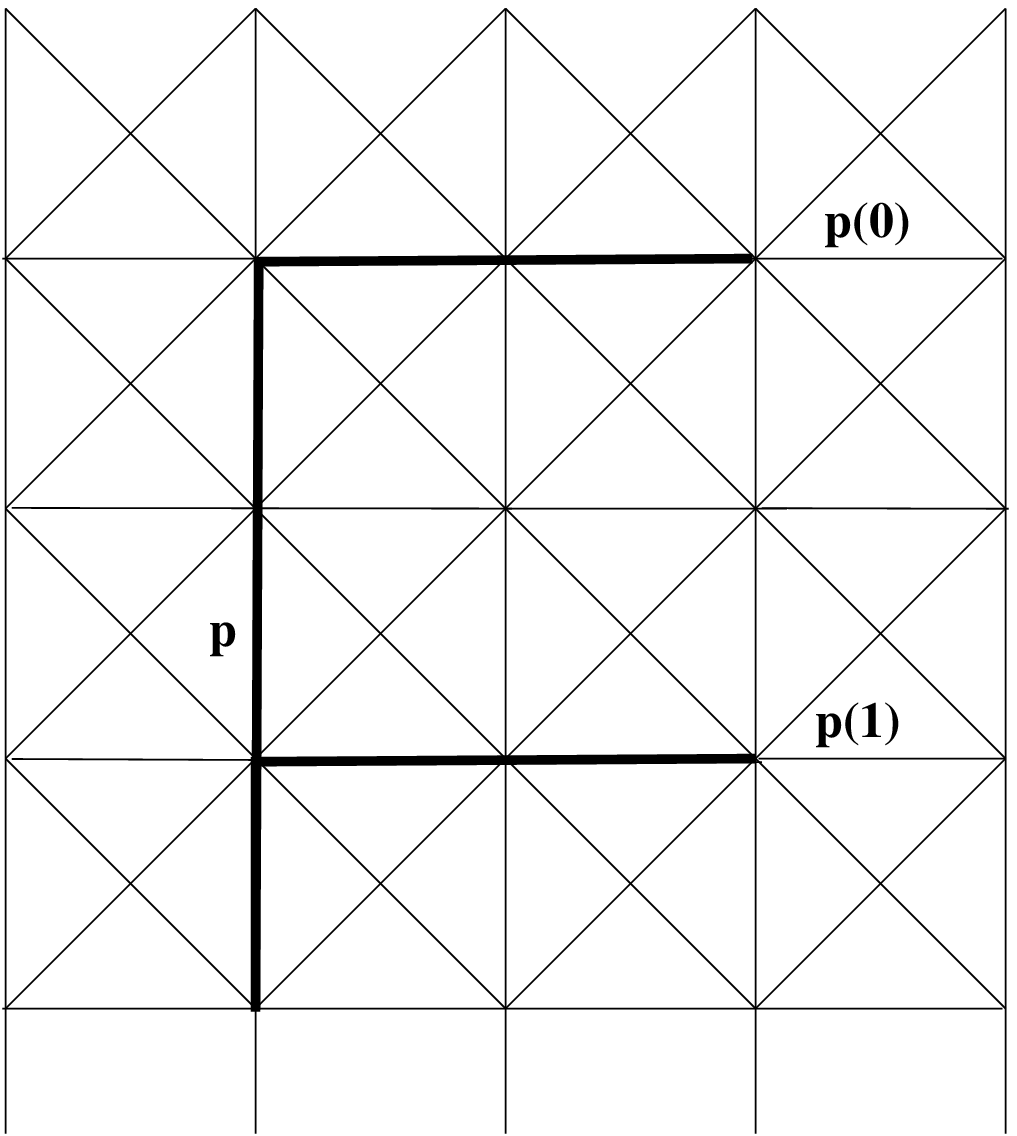}}
\caption{\sl }
\label{type1.fig}
\end{figure}

\begin{figure}[tbh]
\centerline{\epsfxsize=4.5in \epsfbox{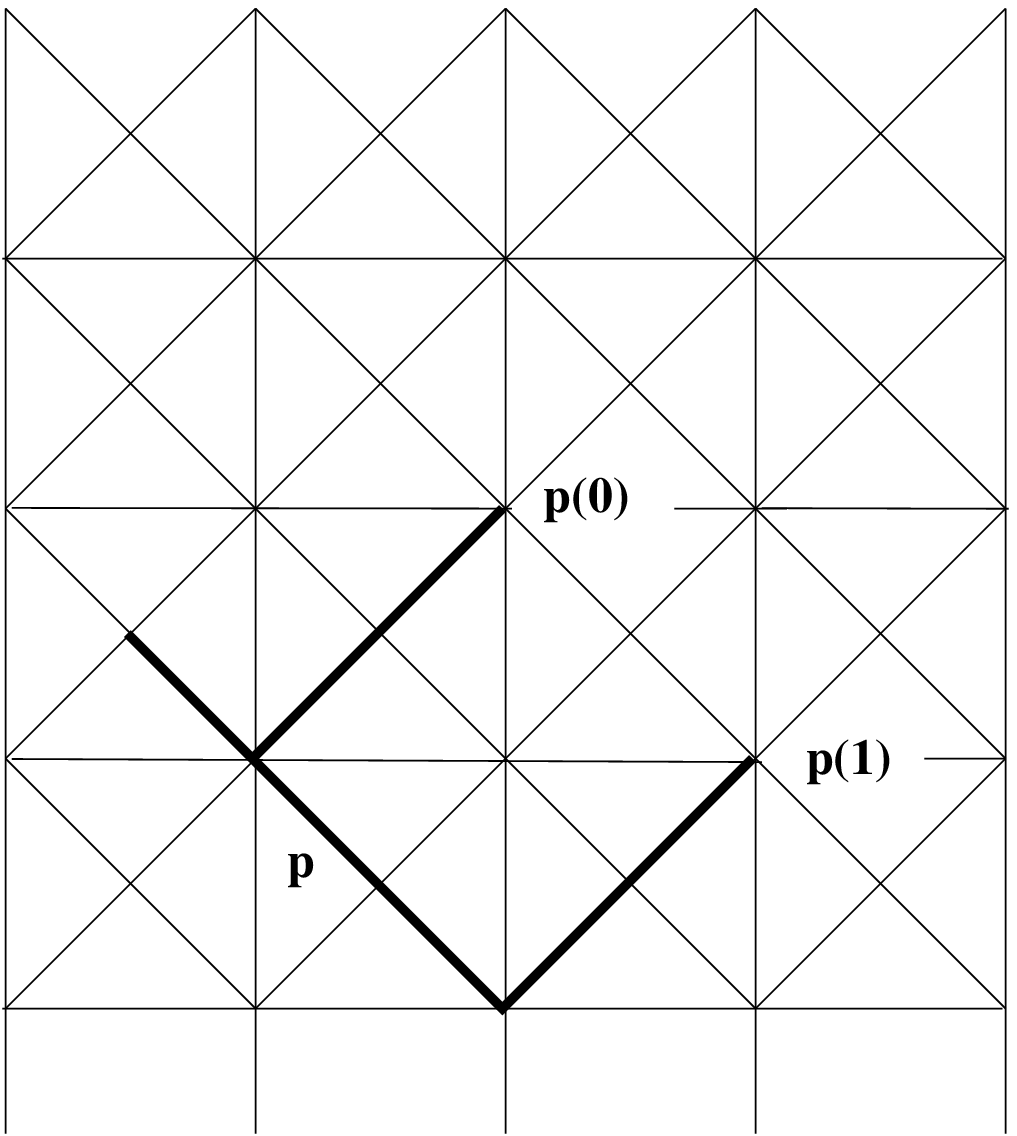}}
\caption{\sl }
\label{type2.fig}
\end{figure}

Accordingly, each $LS(R)$ path $p$ with $\Length(p)\in \N \varpi_1$ has
the shape as in Figure \ref{type1.fig} and each $LS(R)$ path with
$\Length(p)\in \N \varpi_2$ has the shape as in Figure
\ref{type2.fig}. (Some of the geodesic segments in $p$ described
in these figures could have zero length.) 

\begin{lem}
\label{L1}
 Suppose $p$ is an LS path with respect to
$2R$ so that $p(0)$ is a special vertex of $(A, W_{aff})$,
$\Length(p)=\mu\in \Z_+\cdot \varpi_1$ and all breaks of $p$
are at the verticies of $(A, W_{aff})$. Then:

1.  $p$ has breaks only at the special verticies of $(A, W_{aff})$

2. $p_1$ is an LS path with respect to $R$. 

3. If $\mu\in 2\Z_+\varpi_1$, then $p_1(0)- p_1(1)\in Q(R)$.
\end{lem}
\proof 1. Observe that our assumptions imply that the path $p_1$ has
edges parallel to the coordinate axes and it is entirely contained
in the 1-skeleton of the simplicial complex $(A, W_{aff})$. Let
$v$ be a non-special vertex of $(A, W_{aff})$. Then $v$ is
disjoint from all the edges of $(A, W_{aff})$ parallel to the
coordinate axes. Therefore the path $p_1$ is disjoint from the set
of nonspecial verticies. 

2. Suppose that $v=p(t)$ is a break-point of $p$. 
Since $v$ is a special vertex of $(A, W_{aff})$, the chain condition in the definition 
of an LS path (with respect to $R$) at $v$ follows from the chain condition with respect to $2R$. 

3. Since $p\in LS(R)$, it follows from Lemma \ref{trivial} that
$$
p_1(0)- p_1(1)-\mu \in Q(R). 
$$
Since $\mu\in 2\Z_+\cdot \varpi_1\subset Q(R)$, the last assertion of Lemma follows. \qed

\begin{cor}
\label{p_1}
Suppose that $p$ is an $LS(R)$--path with $\Length(p)\in \N\varpi_1$. Then $p$ has
breaks only at the special verticies of $(A, W_{aff})$.
\end{cor}

\begin{lem}
\label{L0} Suppose that $p=p_1*p_2\in LS_1(R)$. Then $p$ can have
at most one break-point which is not a special vertex of $(A,
W_{aff})$; such a point is a point of backtracking of the sub-path
$p_2$.
\end{lem}
\proof By the previous lemma, the path $p_1$ can have breaks only at the special
verticies of $(A, W_{aff})$. Consider the path $p_2$. Suppose that a
break-point $p(t)$ is a {\em turning point} of $p$; let $\tau\in
W_{aff}$ denote the reflection fixing $p(t)$ which sends $p'_-(t)$
to $p'_+(t)$. The vectors $p'_-(t)$, $p'_+(t)\in T_{p(t)}A$ are
tangent to the walls $H_-, H_+$ of $(A, W_{aff})$ which are
parallel to the lines $\{x=y\}, \{x=-y\}$. The linear part of
$\tau$ permutes $H_-, H_+$, hence $\tau$ is the reflection in a
wall which is either vertical or horizontal. Therefore $p(t)$ is a
special vertex of $(A, W_{aff})$.

Hence, if a break-point $p$ is non-special, then the path $p$ has
to backtrack at this point. On the other  hand, our analysis of the shapes of LS paths 
shows that there could be at most one point where $p_2$
backtracks. \qed

\medskip
We now analyze the points of backtracking.

\begin{lem}
\label{special} Suppose that $p=p_2: [0,1]\to V$ belongs to $LS(R)$, $\Length(p)\in \N\varpi_2$, 
$p(t_1)$ is a backtracking point of $p$ which belongs to the wall $\{x=y\}$, and 
$p(0)$ or $p(1)$ belongs to $Q(R)$. Then $p(t_1)$ is a special vertex of $(A, W_{aff})$.
\end{lem}
\proof Suppose that the break-point $p(t_1)=(x_1,x_1)$ is not
special. Consider the line
$$
\ell=\{ x+y=2x_1\}.
$$
Observe that the line $\ell$ contains no verticies of $Q(R)$ since $2x_1$ is an odd number,
while for each $(x, y)\in Q(R)$, the sum $x+y$ is even. Suppose that $p(0)\in Q(R)$. Then
$p(0)\notin \ell$ and hence there exists a geodesic subsegment $\ol{p(0) p(t_2)}$ in $p$ which
is orthogonal to $\ell$ so that $p(t_2)=(x_2,y_2)\in \ell$ is a turning point of $p$. However such segment
clearly cannot contain points in $Q(R)$ since its points $(x,y)$ satisfy the equation
$$
x-y=x_2-y_2\notin 2\Z.
$$
Contradiction. \qed

As a corollary we obtain:

\begin{cor}
\label{C2}
Suppose  that $\si=(\la,\mu,\nu)\in Q(R)^3\cap D_3$ and $\la+\mu+\nu\in 2P(R)$.
Then there exists a path $p\in LS_1(R^\vee)$ contained in $\De$
such that:

1. $p(0)=\la, p(1)=\nu$.

2. All break-points of $p$ are verticies of $(A, W_{aff})$.
\end{cor}
\proof Note that $Q(R)=2P(R^\vee)$. We now apply the results established above to the root
system $R^\vee$; note that this interchanges the roles of short and long fundamental weights,
e.g. the wall $\{x=y\}$ contains the short fundamental weight of $R^\vee$.

Then,  according to the saturation theorem \ref{B2saturation}, $\si \in Tens(R^\vee)$.
Moreover, there exists a path $p\subset \De$ connecting
$\la$ to $\nu$ such that:

1. $\length(p)=(\mu_1, \mu_2)\in 2P(R^\vee) \times 2P(R^\vee)$.

2. $p\in LS_1(R^\vee)$ is a generalized LS path with
respect to $R^\vee$.

Accordingly, the path $p$ is the  concatenation $p=p_1*p_2$ of two
$LS(R^\vee)$--paths $p_i: [0,1]\to V$. Consider the break-points of the path $p$.
According to Corollary \ref{p_1}, all break-points of $p_2$ are special verticies of $(A, W_{aff}^\vee)$.
All break-points of $p_1$ are special verticies of
$(A, W_{aff}^\vee)$ except possibly for a point $p_2(t_1)$ of backtracking which occurs on a wall of $\De$,
where the path $p_1$ intersects $\D \De$ orthogonally (see Theorem \ref{B2saturation} and
Lemma \ref{L0}).
On the other hand, the end-point $\la=p_1(0)$ of the path $p_1$ belongs to $2P(R^\vee)\subset Q(R^\vee)$.
Therefore, it follows from Lemma \ref{special} (applied to the root system $R^\vee$ rather than $R$)
that $p_1(t_1)$ has to be a special vertex of $(A, W_{aff}^\vee)$.

Thus all break-points of $p$ occur in special verticies of $(A, W_{aff}^\vee)$. Since
$P(R^\vee)$ is the set of verticies of $(A,W_{aff})$, all break-points of $p$ are verticies of
$(A, W_{aff})$. \qed

\subsection{Analysis of the $LS(R^\vee)$ paths with singular $\De$-length}
\label{singular}

In Corollary \ref{C2} we have established  that for a large class of triples $(\la,\mu,\nu)\in D_3$,
there exists an $LS_1(R^\vee)$ path $p$ with $\length(p)=\mu$ connecting
$\la$ to $\nu$, so that all break-points of $p$
are verticies of $(A, W_{aff})$.

We now analyze LS paths with respect to $R^\vee$ whose break-points
are verticies of $(A, W_{aff})$. Such paths necessarily belong to $LS(2R)$.
Throughout this section we assume that $p: [0,1]\to V$ is a PL path in
$\De$ so that:

\begin{itemize}
\item $\Length(p)=\mu\in \N\varpi_2$. 

\item $\la:=p(0), \nu:=p(1)\in P(R)$, $\la-\nu\in Q(R)$. 

\item All break-points of $p$ are verticies of $(A,W_{aff})$. 

\item $p\in LS(R^\vee)$. 
\end{itemize}

Unlike the paths $p$ with $\Length(p)\in \N\varpi_1$, the above paths $p$ do not have to belong 
to $LS(R)$.  However, after analyzing these paths, we show that in ``most cases'' they can be replaced 
with a new $LS(R)$ path $\hat{p}$ while keeping the same end-points and the $\De$-length. 

\medskip 
If the path $p$ is not an LS path with respect to $R$ then it has
a break-point which is not a special vertex of $(A,W_{aff})$.

\begin{defn}
Suppose that $v=p(t)$ is a break-point of a PL path $p$. We call
the break-point $x$ {\em illegal} if it violates the axioms of an
LS path (with respect to $R$).
\end{defn}

%{\mini \begin{rem} Observe that at each backtracking point the
%path $p$ satisfies {\em Hecke path condition}, while at each
%turning break point --- it does not.
%\end{rem}}

 Accordingly, we will refer to an illegal break-point which is a
turning/backtracking point as an {\em illegal turn/backtrack}.

\begin{figure}[tbh]
\centerline{\epsfxsize=4.5in \epsfbox{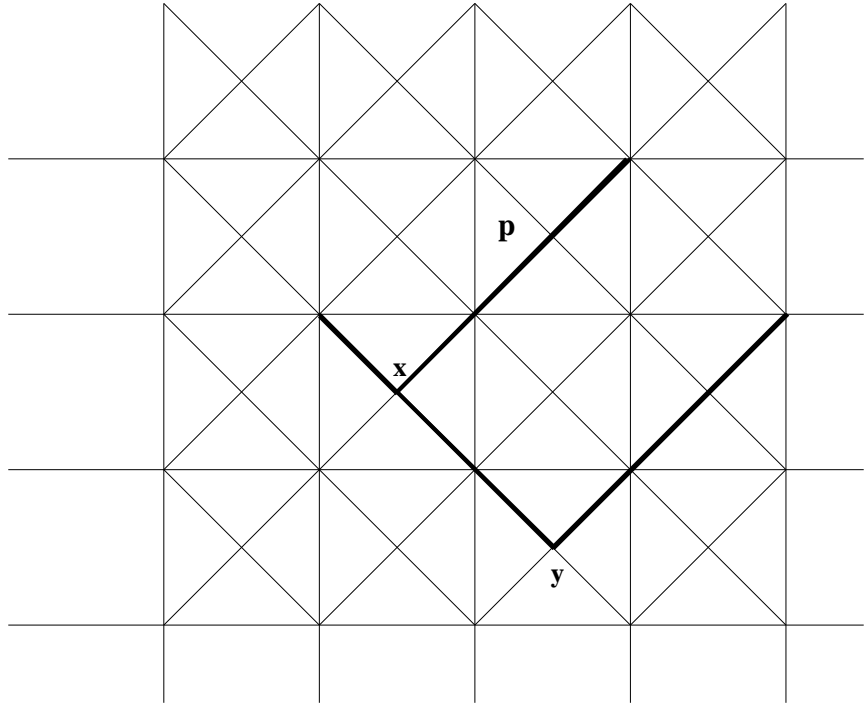}} \caption{\sl
Points $x, y$ are illegal breaks.} \label{illegal.fig}
\end{figure}

\begin{lem}
\label{two}
1. Either $p$ has no illegal breaks or it has 2 illegal turns at
$t_1< t_2$ or one illegal backtrack.

2. In the case of two illegal turns at $t_1, t_2\in [0,1]$, up to
the reversal of orientation, the path $p$ has the shape described
in the Figure \ref{illegal.fig}.
\end{lem}
\proof We first consider the backtracks of $p$. If $p$ has a backtrack at $t$ and $p'_+(t)\notin \De$ 
then 
$$
w(p'_-(t))=p'_+(t) 
$$
where $w\in W_{aff}^\vee$ is a reflection fixing the vertex $p(t)$ and whose linear part is $\tau_1$. 
Since $\tau_1$ is a simple reflection, it follows that the chain 
$$
(p'_-(t), p'_+(t))
$$ 
is necessarily maximal (see Lemma 3.15 in \cite{KM}) and therefore the backtrack at $t$ is 
{\em legal}. Moreover, the reflection $w$ also belongs to $W_{aff}$. 
Thus, if $p$ has an illegal backtrack at a point $t$ then $p'_-(t)\in
-\De, p'_+(t)\in \De$. Therefore $p$ has no other breaks in this
case and we are done. Hence we assume that $p$ has no illegal backtracks.

Let $t_i\in [0,1], i=1,...,m$ denote the illegal turns of
$p$. We set $t_0^+:= 0$, $t_{m+1}^-:= 1$. For each $t_i$ define
$[t_i^-, t_i^+]\subset [0,1]$ to be the smallest subinterval
containing $t_i$ such that $p(t_i^-), p(t_i^+)$ are special
verticies of $(A,W_{aff})$. Then the assumption that for each $i$ the point $p(t_i)$
is an illegal turn  implies that
\begin{equation}
\label{notin}
 p(t_i^+)-p(t_i^-)\notin Q(R).
\end{equation}
On the other hand, the restriction of $p$ to each subinterval $[t_i^+, t_{i+1}^-]$ is
an LS path, therefore
$$
p(t_i^+)- p(t_{i+1}^-)\in Q(R)
$$
(see Lemma \ref{trivial}). Thus
$$
p(0)- p(t_1^-) + p(t_1^+) - p(t_2^-)+... -p(t^-_m)+ p(t^+_m)-
p(1)\in Q(R)
$$
and since $p(0)-p(1)\in Q(R)$, the latter is  equivalent to:
$$
(- p(t_1^-) + p(t_1^+)) + ... +(-p(t^-_m)+ p(t^+_m))\in Q(R).
$$
Thus, since $Q(R)$ has index $2$ in $P(R)$, it follows from
(\ref{notin}) that $p$ has to have an even number of illegal turns.
Since the length of the longest chain in $W\cdot
\varpi_2$ is $4$, the number of illegal turns is $\le 3$, hence it is
either 0 or 2. This proves the first assertion of Lemma.

Suppose that $p$ has illegal turns at $t_1$ and $t_2$, where $t_1<t_2$. Then, analogously
to the case of an illegal backtrack,
$$
p_-'(t_1)\in -\De, \quad p_+'(t_2)\in \De.
$$
Therefore $p\restr_{[0, t_1]}, p\restr_{[t_2,1]}$ are geodesic paths. Moreover,
since the length of the longest chain in $W\cdot \varpi_2$ is $4$,
the path $p$ can have at most one (necessarily legal) break-point
on the open interval $(t_1, t_2)$. Therefore the path $p$ has the
shape as in Figure \ref{illegal.fig}. \qed

\medskip 
Let $p$ be a path as above. Let $x_1=p(t_1),
x_2=p(t_2)$ denote the first and the last illegal breaks of $p$ (possibly $t_1=t_2$). Then 
both breaks occur at nonspecial verticies. According to the above lemma, either $t_1=t_2$ 
and $p$ has an illegal backtrack at this point or both breaks are illegal turns.  
Let $t_0<t_1$ and $t_3>t_2$ be the largest and smallest values of $t$ such
that $p(t)\in P(R)$. Our goal is to show that, with one exception, 
one can always modify the path $p$ on the interval $(t_0, t_3)$ (preserving its $\De$-length)
so that the new path is an $LS(R)$--path.

\begin{figure}[tbh]
\centerline{\epsfxsize=5.7in \epsfbox{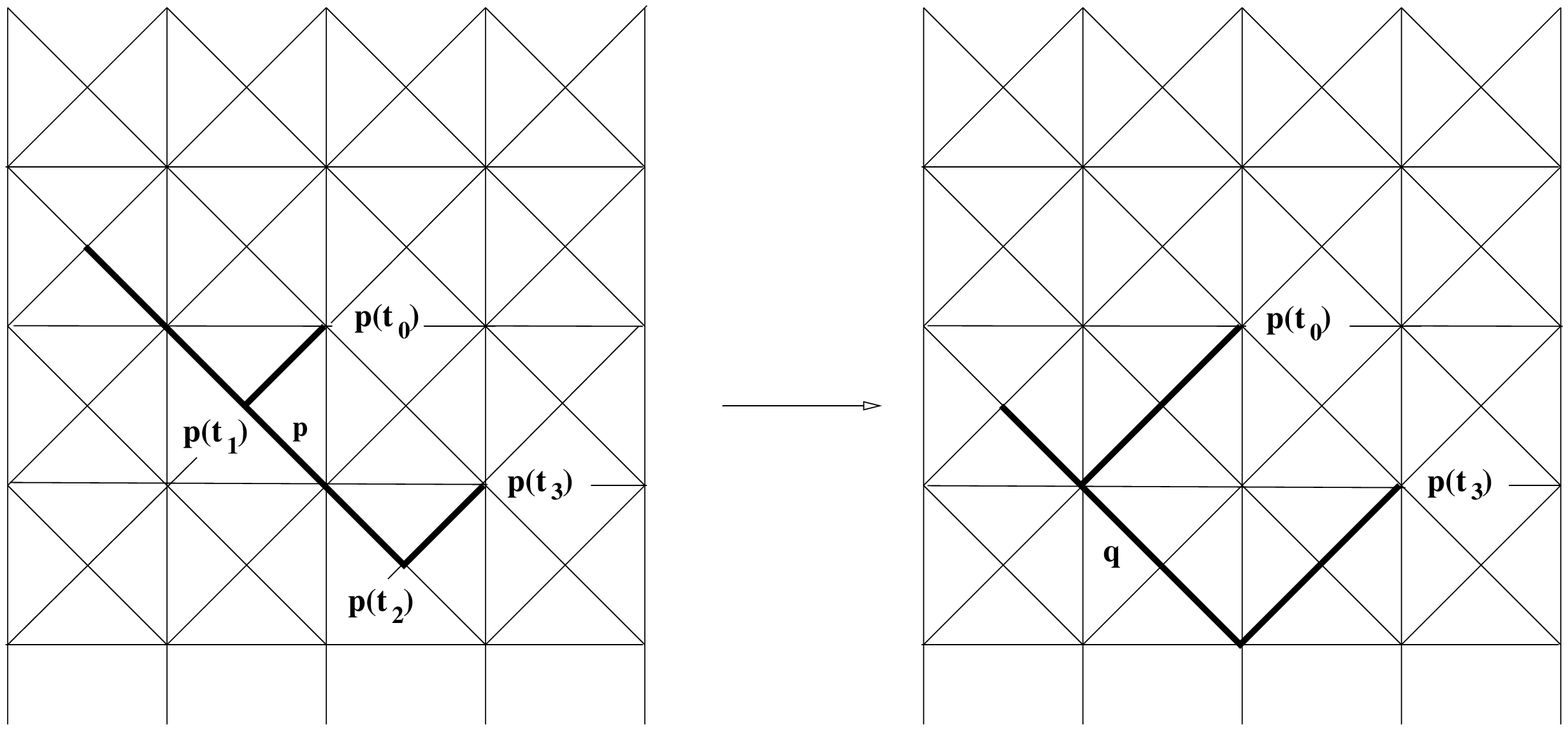}}
\caption{\sl Points $p(t_1), p(t_2)$ are illegal breaks of $p$.}
\label{mod1.fig}
\end{figure}

\begin{lem}
\label{K1}
Suppose that $p(0)\notin \N \varpi_2$ and the image of $p$ is contained in $\De$.
Then there exists a path $q: [t_0, t_3]\to V$ such that:

1.  The concatenation
$$
\hat{p}:=p\restr_{[0, t_0]}* q * p\restr_{[t_3, 1]}
$$
is an $LS(R)$ path contained in $\De$.

2. $\Length(\hat{p})=\Length(p)$.
\end{lem}
\proof There are two cases to consider:

Case 1. The path $p$ restricted to the open interval $(t_1, t_2)$ is not geodesic. In particular, $t_1\ne t_2$ 
and both illegal breaks are illegal turns. 
Then we use the modification described in Figure \ref{mod1.fig}.
It is clear that the new path $\hat{p}$ always satisfies the requirements of Lemma.

\begin{figure}[tbh]
\centerline{\epsfxsize=5.7in \epsfbox{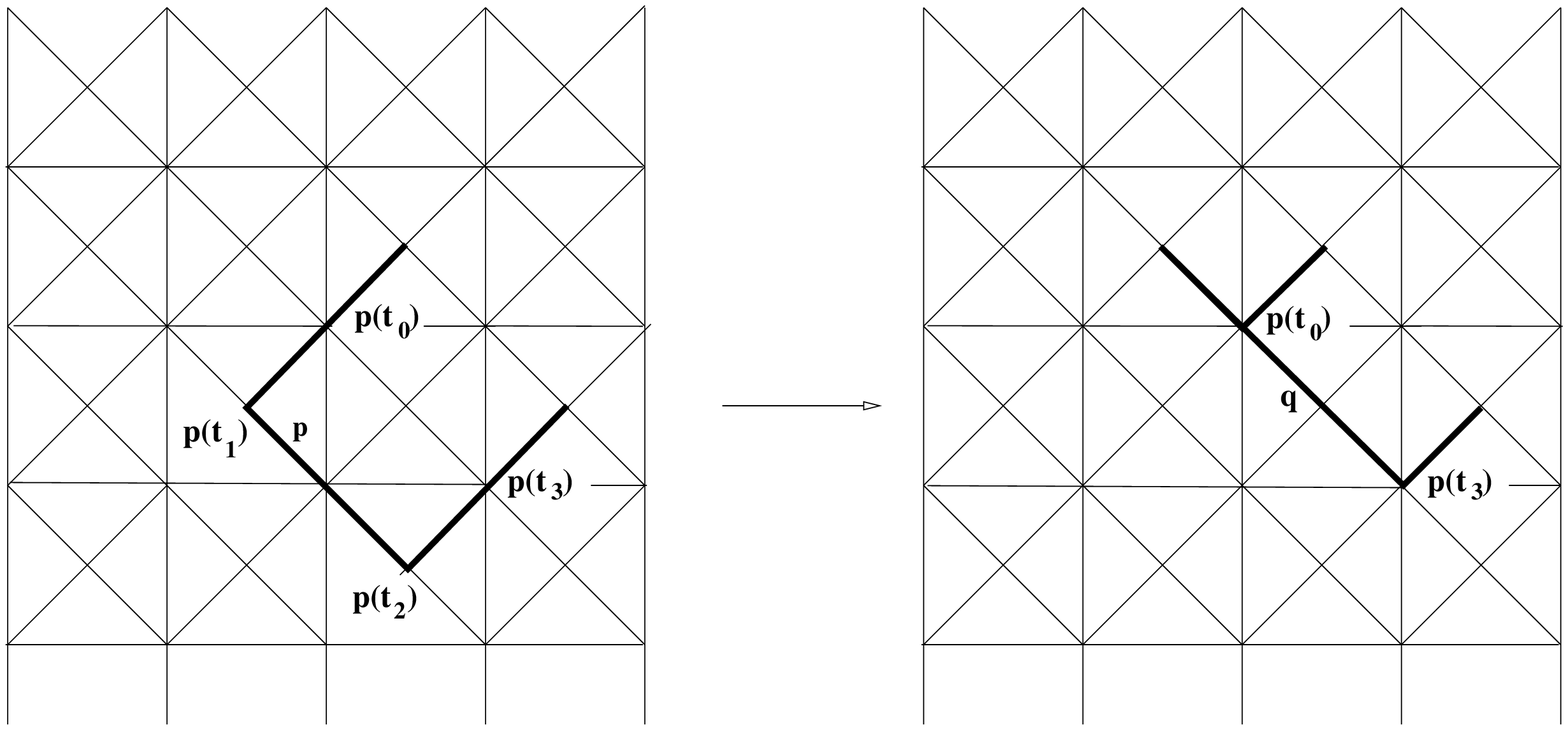}}
\caption{\sl Points $p(t_1), p(t_2)$ are illegal breaks of $p$.}
\label{mod2.fig}
\end{figure}

Case 2. The path $p\restr_{(t_1, t_2)}$ is a (possibly constant) geodesic. Then we use the
modification described in Figure \ref{mod2.fig} by introducing an
extra break $\tau$ between $t_1$ and $t_2$.

It is clear that the new path $\hat{p}$ is an LS path, it has the
same $\De$-length as $p$. However  $\hat{p}$ is not necessarily
contained in $\De$: The point $q(\tau)$ could be outside of $\De$.
This happens if and only if the point $p(t_1)$ lies on the wall
$\{x=y\}$ of $\De$. In this case however $p(0)\in \{x=y\}$ as well which contradicts our hypothesis.
\qed

We now analyze in detail the exceptional case when $\la:=p(0)$ lies on
the wall $\{x=y\}$, i.e. $\la$ belongs
to $\N\varpi_2$. Then the path $p$ has exactly two turns or one
backtrack and all breaks are illegal.

\begin{lem}
\label{K2}
Under the above conditions we have:
$$
\la+\mu+\nu\notin 2P(R).
$$
\end{lem}
\proof Let $p$ have breaks at $t_1\le t_2$. Set $u:= p(t_1), v:=
p(t_2)$. Then
$$
u=\la- (l+\frac{1}{2})\varpi_2, v= \nu- (m+\frac{1}{2})\varpi_2,
$$
$$
\mu=\tau_2((l+\frac{1}{2})\varpi_2 + (v-u) +
(m+\frac{1}{2})\varpi_2),
$$
where $l, m\in \Z$ and $\tau_2(x,y)=(x,-y)$. Note that
$\mu-\tau_2(\mu)\in 2P(R)$. Therefore, modulo $2P(R)$ we have:
$$
\la+\nu+\mu= \la+\nu+(l+\frac{1}{2})\varpi_2 + (v-u) +
(m+\frac{1}{2})\varpi_2= 2\nu+ (2l+1)\varpi_2.
$$
Hence $\la+\mu+\nu$ does not belong to $2P(R)$. \qed

\begin{figure}[tbh]
\centerline{\epsfxsize=5in \epsfbox{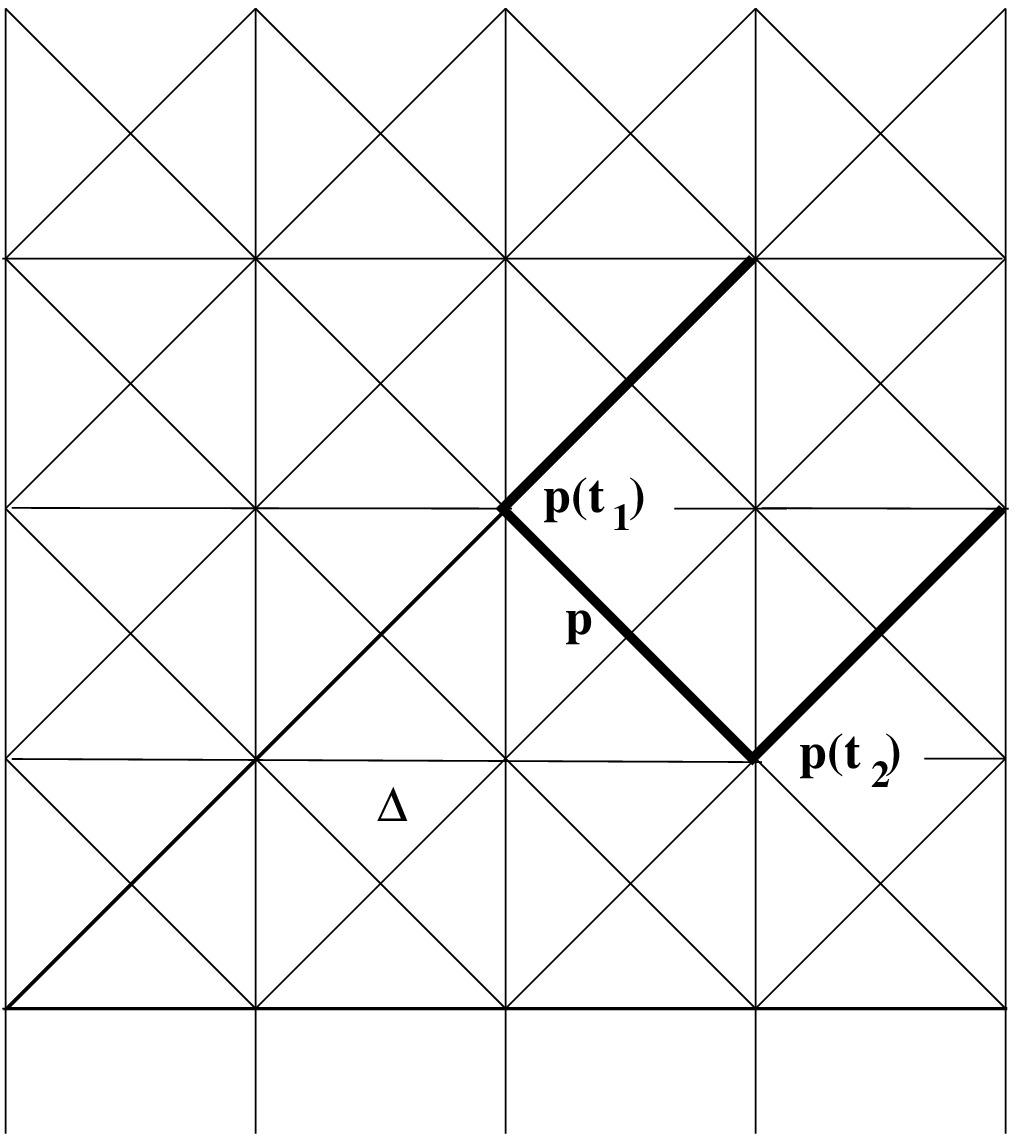}}
\caption{\sl }
\label{case0.fig}
\end{figure}

\begin{lem}
If $\la=n\varpi_2$ and $(\la,\mu,\nu)\in Tens(R)$ then $\la+\mu+\nu\in 2P(R)$. 
\end{lem}
\proof Suppose that $(\la,\mu,\nu)\in Tens(R)$. Then there exists
an LS path $q$ with $\Length(q)=\mu$, $q(0)=\la, q(1)=\nu$, so
that $q$ is entirely contained in $\De$. This path is either
geodesic or has one or two turns, or one backtrack. See Figure
\ref{case0.fig}. We consider the ``generic case'' when $q$ has two turns at the 
points $u=q(t_1), v=q(t_2)$, $t_1<t_2$. Then, analogously to the proof of the previous lemma, 
$$
\mu= (l+s+m)\varpi_2, 
$$
$$
\nu= \la -l\varpi_2 + s\tau_2(\varpi_2) + m\varpi_2,
$$ 
where $u=\la-l\varpi_2, v=\nu-m\varpi_2$ and $l, s, m\in \N$. 
Since $\varpi_2\in P(R)$ and $\tau_2(\varpi_2)-\varpi_2\in 2P(R)$, we obtain:
$$
\nu\equiv \la +l\varpi_2 + s\varpi_2 + m\varpi_2 ~~(\hbox{~mod~} 2P(R)), 
$$
i.e.
$$
\nu \equiv \la+\mu ~~(\hbox{~mod~} 2P(R)). 
$$
Therefore
$$
\la+\mu+\nu\equiv 2\nu\equiv 0 ~~(\hbox{~mod~} 2P(R)). \qed 
$$

We summarize the above results in the following:

%therefore obtain the following:

\begin{prop}
\label{P1}
Suppose that $\mu\in \N\varpi_2$ and $p$ is an LS path with respect to 
$2R$ such that $p$ is contained in $\De$,
$p(0)=\la, p(1)=\nu\in P(R)$, $\Length(p)=\mu$ and all breaks of $p$ 
are at verticies of $(A, W_{aff})$.
Then:

1. $(\la,\mu,\nu)\notin Tens(R)$ if and only if  either
$\la\in \N\varpi_2$ or $\nu\in \N\varpi_2$ and
$$
\la+\mu+\nu\notin 2P(R).
$$

2. Unless $p'(0)\in -\Del, p'(1)\in \Del$, the path $p$ is also an
LS path  with respect to the root system $R$.

3. Unless  $\la$ or $\nu\in \N\varpi_2$, and $\la+\mu+\nu\notin
2P(R)$, there exists a path $\hat{p}$ contained in $\De$ (of
$\De$--length $\mu$) which is an LS path with respect to the root
system $R$, so that $\hat{p}(0)=p(0), \hat{p}(1)=p(1)$.
%whose germs at $0$ and $1$ are equal to those to $p$.
\end{prop}

\subsection{Analysis of $LS_1(R^\vee)$ paths}
\label{regular}

%We now analyze the paths $p$ which are concatenations of
%nontrivial paths $p_1*p_2$ where $\Length(p_i)\in \N\varpi_i$. We
%assume that:

%1. $p$ belongs to $LS_1(2R)$.

%2. All breaks of $p$ are at verticies of $(A, W_{aff})$.

%3. $p(0)=\la, p(1)=\nu$ are both not on the wall $\{x=y\}$.

%and $\la, \nu\notin \N\varpi_2$.

In the previous section we proved that for ``most'' singular $LS(2R)$ paths $p_i\subset \De$ (with 
$\Length(p_i)\in \N\varpi_i$), whose break-points 
are verticies of $(A, W_{aff})$, we can replace  $p_i$ with a new path $\hat{p}_i$ which has the same $\De$-length, 
same end-points and is still contained in $\De$ (Proposition \ref{P1}). 
The goal of this section is to prove a similar statement 
for paths $p=p_1*p_2\in LS_1(2R)$. The naive idea would be to replace each $p_i$ with $\hat{p}_i$ 
using Proposition \ref{P1} and then take $\hat{p}:= \hat{p}_1* \hat{p}_2$. The are two issues 
however which have to be addressed:

(1) It might happen that the path $p_2$ is ``exceptional'' from the point of view of Proposition \ref{P2}, 
i.e. $p_2(0)\in \{x=y\}$. 

(2) We have to ensure that at the concatenation point between $\hat{p}_1$ and $\hat{p}_2$ 
the new path satisfies the axiom of an $LS_1(R)$ path. 

\begin{figure}[tbh]
\centerline{\epsfxsize=5.7in \epsfbox{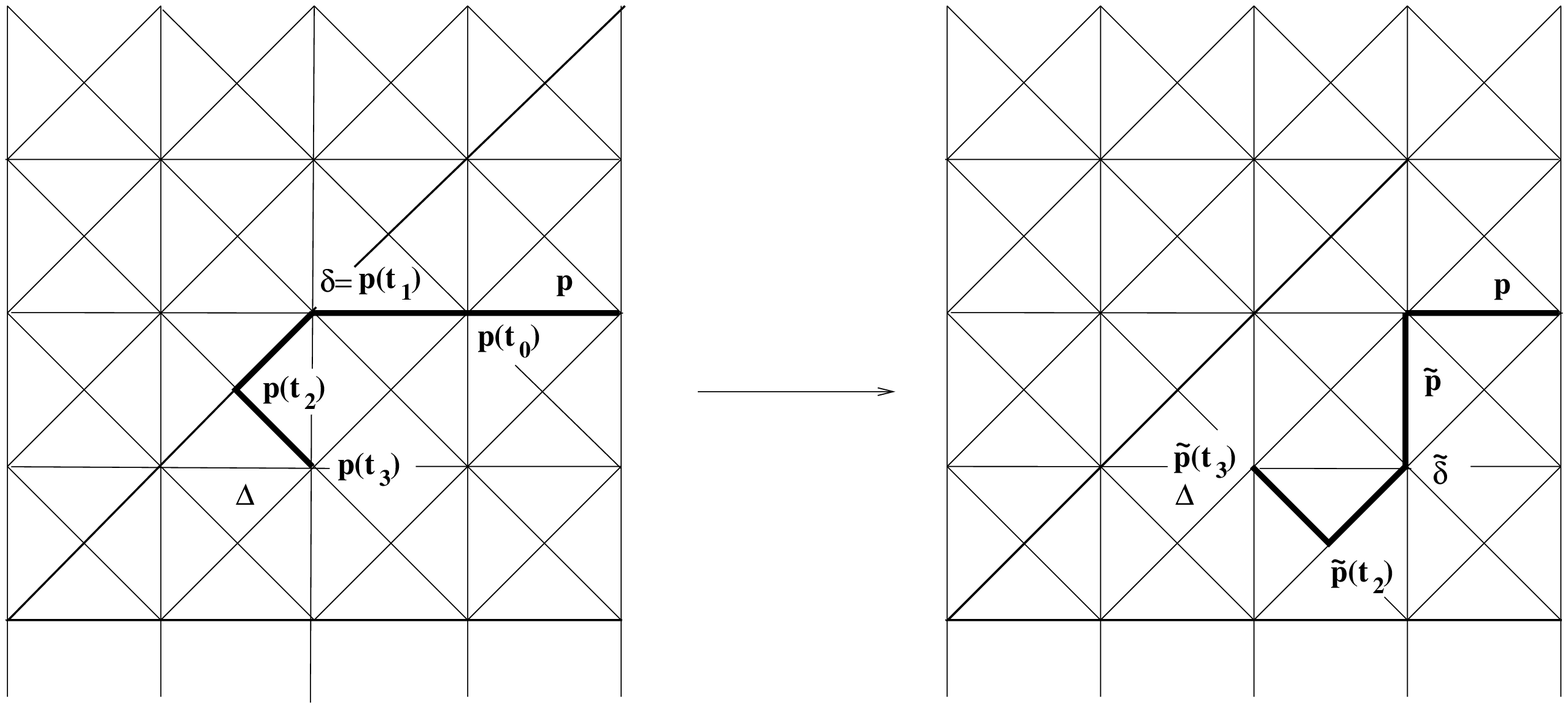}}
\caption{\sl }
\label{mod3.fig}
\end{figure}

It turns out that the issue (1) is trickier to handle: We cannot use Proposition \ref{P1} directly and 
are forced first to change the `concatenation point'' (Figure \ref{mod3.fig}) and move it away from 
the wall $\{x=y\}$.

The main result of this section is the following proposition: 

\begin{prop}
\label{P2} Suppose that $\si=(\la, \mu, \nu)\in P(R)^3$ is such
that:

1. $\la, \nu\notin \N\varpi_2$.

2. There exists a path $p\in LS_1(2R)$ from $\la$ to $\nu$, which is contained in $\De$,
all whose break-points are verticies of $(A, W_{aff})$ and   
so that 
$$
\length(p)=(\mu_1,\mu_2), \mu=\mu_1+\mu_2.
$$
Then $\si\in Tens(R)$.
\end{prop}
\proof We start by analyzing the path $p$. Our goal is to replace it with a a new path which is in $LS_1(R)$ 
and which still satisfies condition 2. 

Set $\del:= p_1(1)$. According to Lemma \ref{L1}, the path $p_1$
is an LS path with respect to $R$. In particular, $\del\in P(R)$.
On the other hand, if $p_2'(0)$ does not belong to $-\De$, then,
according  to the second part of Proposition \ref{P1}, the path
$p_2$ is an LS path with respect to $R$. Hence $p=p_1*p_2$ belongs
to $LS_1(R)$, since the vertex $\del$ is special and the
generalized chain condition at this point (with respect to $R$)
follows from the generalized chain condition at this point (with
respect to $2R$).

\medskip
We now consider the case $p_2'(0)\in -\De$. Observe that, since
$\del$ is a special vertex, for each $\eta\in -\De$ and every
$w\in W$,
$$
\eta\ge w(\eta),
$$
it follows that for every $LS(R)$-path $q_2$, the concatenation
$$
p_1*q_2
$$
belongs to $LS_1(R)$.

\medskip
Case 1. $p_1(1)\notin \{x=y\}$. Then, according to Part 3 of
Proposition \ref{P1}, there exists an $LS(R)$--path $q_2$
(entirely contained in $\De$) starting at $\del$, ending at $\nu$,
with $\Length(q_2)=\mu_2$.  Hence the concatenation
$\hat{p}:=p_1*q_2$ is a generalized LS path with respect to the
root system $R$, $\hat{p}$ is contained in $\De$, and therefore
$(\la,\mu,\nu)\in Tens(R)$.

\begin{figure}[tbh]
\centerline{\epsfxsize=5in \epsfbox{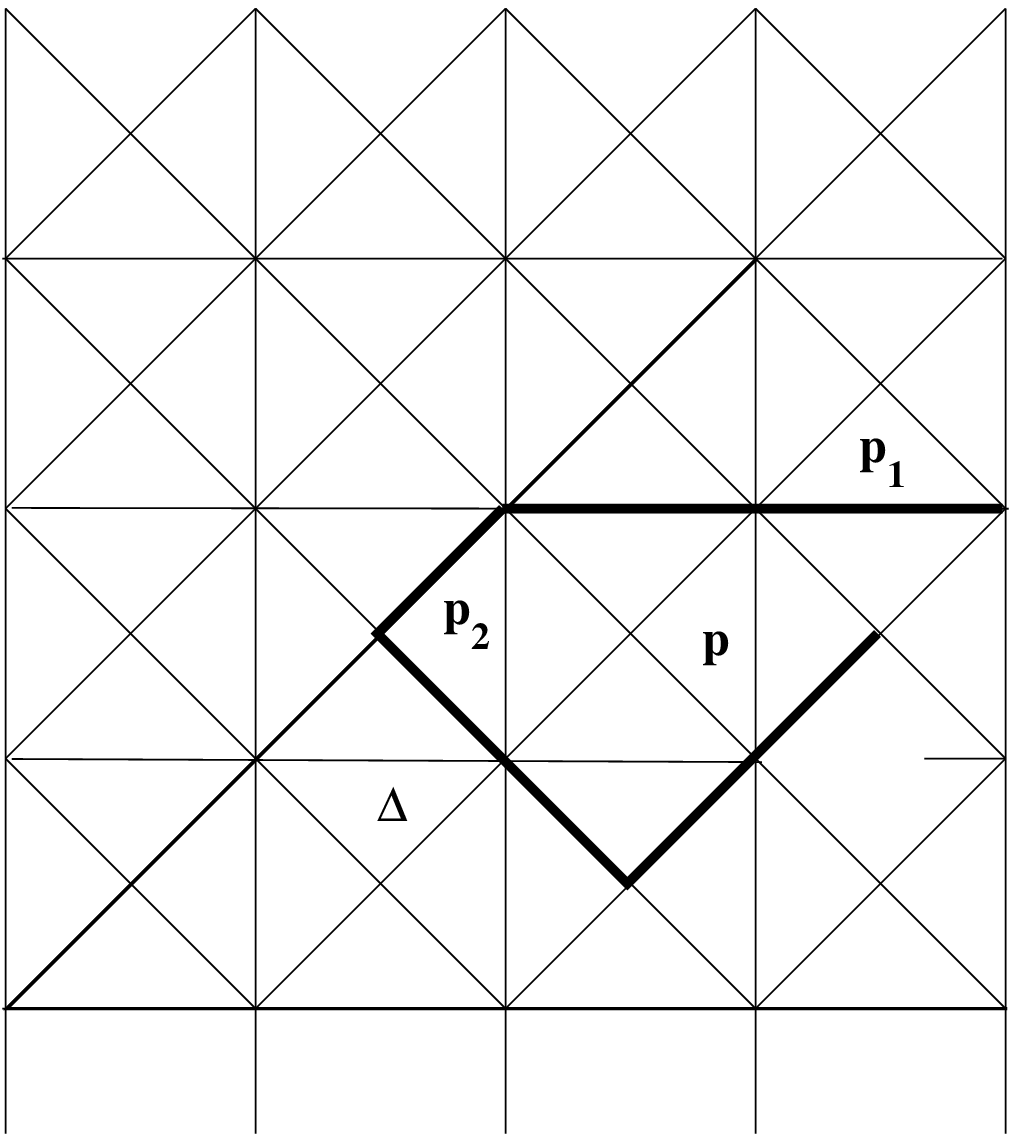}}
\caption{\sl }
\label{shape.fig}
\end{figure}

\medskip
Case 2. $p_1(1)\in \{x=y\}$. Since $p_2'(0)\in -\De$, then
$p_1'(1)\in \R_- \times \R_-$. Since $\del=p_1(1)$ belongs to the
wall $\{x=y\}$ and $p_1$ is contained in $\De$, it follows that
$p_1'(1)\in -\De$. Hence $p_1$ is a geodesic path and the
entire path $p$ has the shape as in Figure \ref{shape.fig}.

We let $t_1\in [0,1]$ be such that $\del=p(t_1)$ is the
concatenation point, let $t_0<t_1$ be the maximal value of $t$
such that $p(t)$ is a special vertex. Let $t_2>t_1$ be the first
value of $t$ where $p(t)$ is not geodesic, $t_3>t_2$ be the first
value of $t$ such that $p(t_2)$ is a (special) vertex. We now
replace the restriction $p\restr_{[t_0, t_3]}$ with the new path
$\t{p}: [t_0, t_3]\to V$ described in Figure \ref{mod3.fig}.  
Observe that $p(t_3)\notin \D\De$, for otherwise 
the path $p_2$ has exactly one illegal turn which contradicts Lemma \ref{two}. 
Therefore the path $\tilde{p}$ is contained in $\De$. 

Moreover, $\t{p}(t_0)=p(t_0), \t{p}(t_3)=p(t_3)$,
$$
\length(\t{p})= \length(p\restr_{[t_0, t_3]}) 
$$
Thus we define the path
$$
q:= p\restr_{[0, t_0]}* \t{p}[t_0, t_3] *  p\restr_{[t_3, 1]}.
$$
Note that the new path has a (legal) turn at the point
$\t{p}(t_3)$ and in addition, two (illegal) simple turns.
Therefore $q=q_1*q_2$  is still not a generalized LS
path. However it has the property that its concatenation point
$\t{\del}=q(t_1)$ is not on the wall $\{x=y\}$. Thus we
have reduced the argument to Case 1 and hence $\si \in Tens(R)$.

This concludes the proof of Proposition \ref{P2}. \qed

\begin{cor}
\label{C3} Suppose that $\si=(\la, \mu, \nu)\in Q(R)^3\cap D_3$ and $\la+\mu+\nu\in 2P(R)$.
Then $\si\in Tens(R)$.
% unless two of the three vectors $\la, \mu, \nu$ belong to
%$\N\varpi_2$ and $\la+\mu+\nu\notin 2P(R)$.
\end{cor}
\proof If all three vectors $\la, \mu, \nu$ do not belong to
$\N\varpi_2$ then the assertion follows immediately from the
combination of Corollary \ref{C2} and Proposition \ref{P2}. If
one of these vectors belongs to $\N\varpi_2$ (by permuting $\la, \mu, \nu$ we can assume that this vector 
is $\mu$), then we use Corollary \ref{C2} and Proposition \ref{P1}. \qed

\subsection{Proof of Theorem \ref{mainBC}}
\label{sect:BC}

Our final goal is to eliminate the assumption that $\la, \mu,
\nu\in Q(R)$ in Corollary \ref{C3}. Recall that $(\la, \mu, \nu)\in
D_3\cap (P(R))^3$, $\la+\mu+\nu\in Q(R)$ and either at most one of the vectors
$\la, \mu, \nu$ belongs to $\N\varpi_2$ or $\la+\mu+\nu\in 2P(R)$. If one of the vectors 
$\la, \mu, \nu$ belongs to $\N \varpi_2$, we can assume (by relabelling) that this vector is 
$\mu$. (Note that $\varpi_2\in Q(R)$ which makes this relabelling  consistent with our convention \ref{convent}.) 

Observe that $\la, \mu, \nu\in P(R)=Q(R^\vee)$, $\la+\mu+\nu\in Q(R)=2P(R^\vee)$.
Therefore, since the root systems $B_2$ and $C_2$ are isomorphic, we can
apply Corollary \ref{C3} to the triple $(\la,\mu,\nu)$ {\em with
respect to the coroot system $R^\vee$} and conclude that
the triple $(\la,\mu,\nu)$ belongs to $Tens(R^\vee)$.

Hence there exists a generalized  LS path $p$ in $\De$ with respect to
the root system $R^\vee$, connecting $\la$ to $\mu$ and having
$$
\length(p)= (\mu_1, \mu_2).
$$
This path is a concatenation $p=p_1*p_2$ of $LS(R^\vee)$ paths $p_1, p_2$, 
the edges of the path $p_1$ are parallel to the $x$ and $y$ axes.

\begin{rem}
Strictly speaking, since in the root system $R^\vee$ the long and the short roots are reversed, 
according to our notation \ref{note}, we would have to use the concatenation 
$p_2*p_1$ rather than $p_1*p_2$. However reversing the roles of $\la$ and $\nu$ 
eliminates this problem.
\end{rem}

\begin{lem}
The breaks of $p$ occur only at special verticies of the Coxeter complex $(A, W_{aff}^\vee)$.
\end{lem}
\proof Observe that our analysis of generalized LS paths (see Lemma
\ref{L0} applied to the root system $R^\vee$) shows that the break-points of $p$ occur only at special
verticies of $(A, W_{aff}^\vee)$ with the sole exception of  a
single break of the sub-path $p_1$ which can occur at a nonspecial
vertex $p(t_1)$, and where $p$ {\em backtracks} and has germ parallel to the $y$ axis.

We claim that this is impossible. Let $[t_0, t_2]$ be the smallest
interval containing $t_1$ such that $p(t_0)=p(t_1)\in P(R^\vee)$.
We then replace the path $p$ with the path $\t{p}$ by {\em
eliminating this backtracking}:
$$
\t{p}= p\restr_{[0,t_0]} * p\restr_{[t_1,1]}.
$$
 Then
$$
\t{\mu}:=\Length(\t{p})=\Length(p)- \varpi_1.
$$
It is clear that the new path $\t{p}$ is a generalized LS path
with respect to the root system $R^\vee$ and moreover its breaks occur only at the special verticies  
of $(A, W_{aff}^\vee)$, i.e. at verticies of $(A, W_{aff})$. We are now in position to apply Proposition 
\ref{P2} (with respect to the root system $R$).  

1. Either $(\la, \t\mu, \nu)\in Tens(R)$,  

2. Or two of the vectors $\la, \mu, \nu$ belong to $\N\varpi_2$. 

In Case 1, since $(\la, \t\mu, \nu)\in Tens(R)$, it follows that
$$
\la+\mu+\nu-\varpi_1=\la+\t{\mu}+\nu\in Q(R).  
$$
However  $\la+\mu+\nu\in Q(R)$ (by the assumption in Theorem \ref{mainBC}) and 
$\varpi_1\notin Q(R)$. Contradiction. 

In Case 2, as it was observed in the beginning of this section, 
we can assume that $\mu\in \N \varpi_2$. Then the path $p$ cannot contain a subsegment 
parallel to the $y$-axis, i.e. the above backtracking in the path $p$ is impossible.  \qed

\begin{cor}
\label{C4}
There exists a path $p$ in $\De$,
connecting $\la$ to $\nu$, which is a generalized LS path with
respect to $R^\vee$ all whose breaks are at verticies of $(A,
W_{aff})$ and such that $\length(p)=(\mu_1,\mu_2)$.
\end{cor}

We now can finish the proof of Theorem \ref{mainBC}. Consider a path 
$p$ in $\De$ as in Corollary \ref{C4}. The breaks in this path occur only 
in verticies of $(A, W_{aff})$. 
Suppose that at most one of the vectors $\la, \mu, \nu$ belongs to 
$\N\varpi_2$. Then, according to Proposition \ref{P2}, $\si\in Tens(R)$. 

Consider the exceptional case, say, $\la, \mu\in \N\varpi_2$. 
Then, according to the hypothesis of Theorem \ref{mainBC}, 
$\la+\mu+\nu\in 2P(R)$. Therefore we can apply Proposition \ref{P1} and $\si\in Tens(R)$. \qed

\medskip Below we express $Tens(R)$ as a union of elementary sets, where $G=Sp(4,\C)$. 
In what follows, $\Z_+=\{0, 1, 2, ...\}$. Let
$$
E_1=\{(\la, \mu, \nu)\in L^3: \la+\mu+\nu\in 2P(R), \la\in \Z_+
\varpi_2, \mu\in \Z_+ \varpi_2\}\cap \P(G),
$$
$$
E_2=\{(\la, \mu, \nu)\in L^3: \la+\mu+\nu\in 2P(R), \la\in \Z_+
\varpi_2, \nu\in \Z_+ \varpi_2\}\cap \P(G),
$$
$$
E_3=\{(\la, \mu, \nu)\in L^3: \la+\mu+\nu\in 2P(R), \nu\in \Z_+
\varpi_2, \mu\in \Z_+ \varpi_2\}\cap \P(G),
$$
$$
E_1'=\{(\la, \mu, \nu)\in \La: \la\notin \Z_+ \varpi_2, \mu\notin \Z_+
\varpi_2\}\cap \P(G),$$
$$
E_2'=\{(\la, \mu, \nu)\in \La: \la\notin \Z_+ \varpi_2, \nu\notin \Z_+
\varpi_2\}\cap \P(G),$$
$$
E_3'=\{(\la, \mu, \nu)\in \La: \mu\notin \Z_+ \varpi_2, \nu\notin \Z_+
\varpi_2\}\cap \P(G).$$ 
Then
$$
Tens= \bigcup_{i=1}^3 E_i \cup \bigcup_{i=1}^3 E'_i.
$$

\section{Computation of $Tens(G_2)$}
\label{G2}

Let $R$ be the root system $G_2$, and let $L:= P(R)$ denote the weight lattice.
We  let $\varpi_1, \varpi_2$ denote the fundamental weights of $R$ so that
 $\varpi_2$ is the longer weight. Let $H_i$ denote the walls $\R\varpi_i$, $i=1, 2$.  
 We will use the coordinates $[x,y]$ for vectors
$\la=x\varpi_1+ y\varpi_2$ in $V=P(R)\otimes \R$, so that the chamber
$\De$ is given by the inequalities $x\ge 0, y\ge 0$. Let $G$ be the complex
simple Lie group with the root system $R$ and maximal compact subgroup $K$.
Recall that $\P(G)=D_3(G/K)\subset \De^3$ denotes the convex cone given by the stability and chamber inequalities.
Note that the permutation group on 3 elements $S_3$ acts on $\P(G)$ (by permuting $\la_1, \la_2, \la_3$)
and this action preserves $Tens(G)$.

The following theorem gives a complete description of the semigroup $Tens(G)$.

\begin{thm}\label{g2}
Suppose that $\si=(\la_1, \la_2, \la_3)\in \P(G)\cap L^3$. Then:

1. If at most one of the vectors $\la_i$ is a multiple of $\varpi_2$ then $\si$
belongs to $Tens(G)$.

2. Suppose that $\la_1=y_1\varpi_2, \la_2=y_2\varpi_2$. Then $\si\notin Tens(G)$ if and only if
$\si$ belongs to the union $\E_1\cup \E_2 \cup \E_3$ of the following ``exceptional'' elementary sets:
$$
\E_1= \left\{ \left(\left[\begin{array}{c}
0\\
y_1\end{array}\right],
\left[\begin{array}{c}
0\\
y_2\end{array}\right],
\left[\begin{array}{c}
1\\
y_3\end{array}\right]\right): y_1, y_2, y_3\in \Z_+\right\},
$$
$$
\E_2= \left\{ \left(\left[\begin{array}{c}
0\\
1+n+m\end{array}\right],
\left[\begin{array}{c}
0\\
1+n+2m\end{array}\right],
\left[\begin{array}{c}
1+3m\\
0\end{array}\right]\right): n, m\in \Z_+\right\},$$
$$
\E_3= \left\{ \left(\left[\begin{array}{c}
0\\
1+n+m\end{array}\right],
\left[\begin{array}{c}
0\\
1+m\end{array}\right],
\left[\begin{array}{c}
1+3m\\
1+n\end{array}\right]\right): n, m\in \Z_+\right\}.
$$
\end{thm}

\noindent %Here and in what follows $\Z_+=\{z\in \Z: z\ge 0\}$.
Note that the sets $\E_2, \E_3$ can be also described as follows.
Let $\phi_i(\si)=\phi_i(\la_1, \la_2, \la_3)$, $i=0, 1, 2$
be given by
$$
\phi_0(\si)=(2x_1-x_2+x_3) +3(y_1-y_2+y_3),
$$
$$
\phi_1(\si)= (x_1+x_2-x_3)+  (y_1+ 2y_2 - y_3),
$$
$$
\phi_2(\si)= (x_1+x_2-x_3)+ 3(y_1+y_2-y_3).
$$
The inequalities $\phi_i(\si)\ge 0, i=0, 1, 2$, appear in the system of stability inequalities defining $\P(G)$
(see \cite{KLM1}). Then
$$
\E_2=\{\si=(\la_1, \la_2, \la_3): \la_1, \la_2\in \Z \varpi_2, \la_3\in \Z \varpi_1, \phi_0(\si)=1\},
$$
$$
\E_3=\{\si=(\la_1, \la_2, \la_3): \la_1, \la_2\in \Z \varpi_2, \phi_i(\si)=1, i=1, 2\}.
$$
Thus $\E_2, \E_3$ are sets of lattice points (i.e. elements of $L$)
in translates of strata of the boundary of the cone $\P(G)$.

We get the following corollaries of the above theorem:

\begin{cor}
\label{c1}
If $\si\in \P(G)\cap L^3$ is a nonsingular triple then $\si\in Tens(G)$.
\end{cor}

\begin{cor}
In the decomposition of $Tens(G)$ as the union of elementary sets, the elementary sets are given by
inequalities only and there are no congruence conditions.
\end{cor}

\proof (of Theorem \ref{g2}). The proof that each triple $\si$ which is not in $S_3\cdot (\E_1\cup \E_2\cup \E_3)$
belongs to $Tens(G)$ is a, rather uninteresting, computation. The proof that $(\E_1\cup \E_2\cup \E_3)\cap
Tens(G)=\emptyset$ is based on  analysis of LS paths (and generalized LS paths) with respect to the root system $G_2$.

\begin{figure}[tbh]
\centerline{\epsfxsize=3in \epsfbox{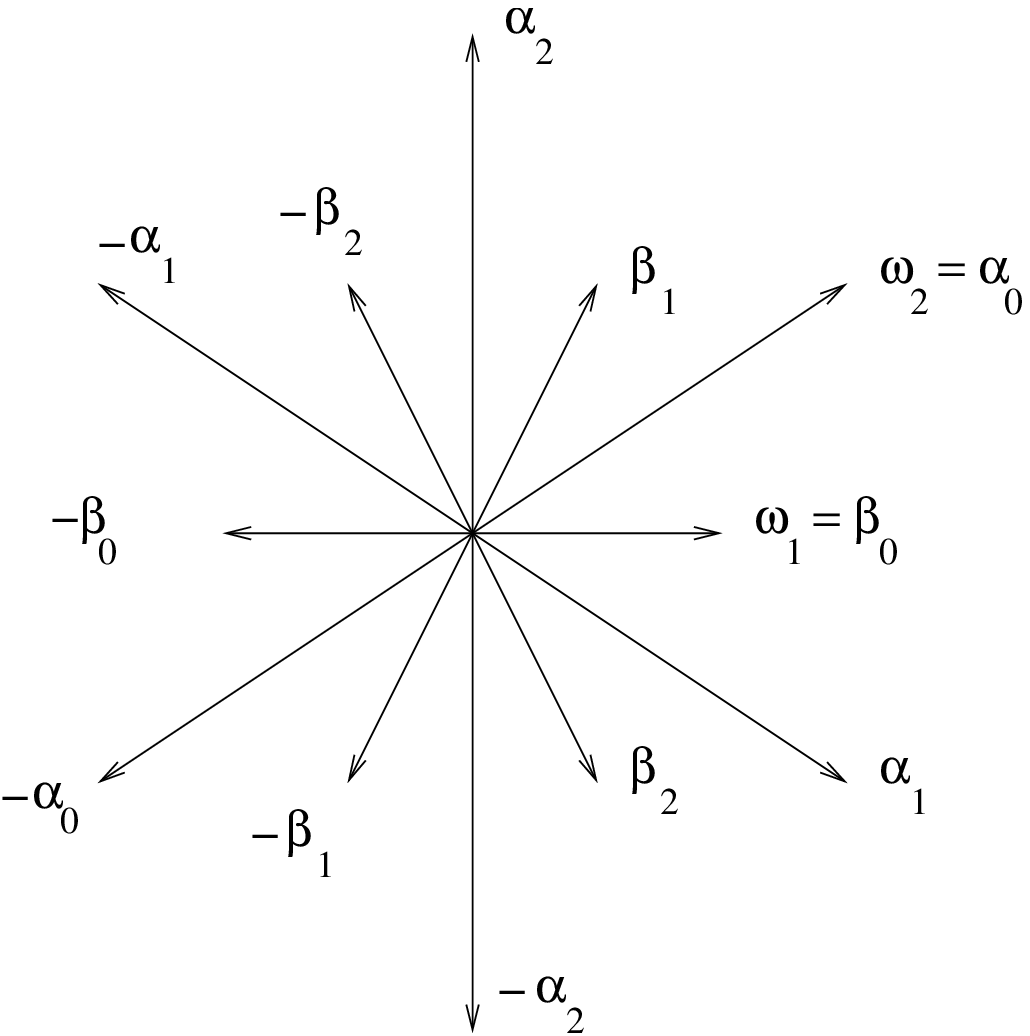}}
\caption{\sl }
\label{f0}
\end{figure}

\begin{notation}
We define the following elements of $L$ (see Figure \ref{f0}):
$$
\al_0:=\varpi_2, \al_1:=3\varpi_1-\varpi_2, \al_2:= \al_0-\al_1.
$$
$$
\be_0:= \varpi_1, \be_1:= \varpi_2-\be_0, \be_2:= \be_0-\be_1.$$
\end{notation}

The following observation will be very useful for the proofs of Propositions \ref{p2}, \ref{p3} below:

\begin{observation}
\label{ob}
Let $(A, W_{aff})$ be the affine Coxeter complex for the root system $R=G_2$ and $W$ be the finite 
Weyl group of $R$. Suppose that $W'$ is the stabilizer of a vertex $v$ in $(A, W_{aff})$; 
we identify $W'$ with a subgroup of $W$. Let $\eta_0\ge \eta_1\ge ...\ge \eta_m=\al_0$ 
is a $W'$-chain which is maximal as a $W$-chain. Then $W'=W$ and $v$ is a special vertex. 
\end{observation}

\begin{defn}
Let $p: [0,1]\to V$ be a Hecke path. We say that $p$ is {\em modeled on a chain}
$\eta_0\ge ...\ge \eta_m$ if the following holds:

Let $t_0=0< ...< t_m<1$, where $t_1,...,t_m$ denote
the break-points of $p$. Since $p$ is a Hecke path, the sequence
$\zeta_0=p'(t_0), \zeta_1=p_+'(t_1), ..., \zeta_m=p'_+(t_m)$
is a chain. We then require each $\zeta_i$ to be a multiple of $\eta_i$, $i=0,...,m$.
\end{defn}

We let $w=(23), u=(13), v=(312)$ denote permutations in $S_3$.

 \begin{prop}\label{p0}
 Suppose that $\si=(\la_1,\la_2,\la_3)\in \P(G)\cap L^3,$
 is such that at most one vector $\la_i$ is a multiple of
 $\varpi_2$. Then $\si\in Tens(G)$.
 \end{prop}
 \proof  Recall that in \cite{KLM3} we have  computed the semigroup generators for $\P(G)\cap L^3$. These are
 the following triples $\del_i, \eps_j, i=1,...,9, j=1, 2$, and their images under the $S_3$--action:

 $$
\del_1:=\left(\left[\begin{array}{c}
1\\
0\end{array}\right],
\left[\begin{array}{c}
1\\
0\end{array}\right],
\left[\begin{array}{c}
0\\
0\end{array}\right]\right),
\quad
\del_2:=\left(\left[\begin{array}{c}
0\\
1\end{array}\right],
\left[\begin{array}{c}
0\\
1\end{array}\right],
\left[\begin{array}{c}
0\\
0\end{array}\right]\right),
$$
$$
\del_3:=\left(\left[\begin{array}{c}
1\\
0\end{array}\right],
\left[\begin{array}{c}
1\\
0\end{array}\right],
\left[\begin{array}{c}
1\\
0\end{array}\right]\right),
\quad
\del_4:=\left(\left[\begin{array}{c}
0\\
1\end{array}\right],
\left[\begin{array}{c}
0\\
1\end{array}\right],
\left[\begin{array}{c}
0\\
1\end{array}\right]\right),
$$
$$
\del_5:=\left(\left[\begin{array}{c}
0\\
1\end{array}\right],
\left[\begin{array}{c}
0\\
1\end{array}\right],
\left[\begin{array}{c}
3\\
0\end{array}\right]\right),
\quad
\del_6:=\left(\left[\begin{array}{c}
0\\
1\end{array}\right],
\left[\begin{array}{c}
0\\
2\end{array}\right],
\left[\begin{array}{c}
3\\
0\end{array}\right]\right),
$$
$$
\del_7:=\left(\left[\begin{array}{c}
0\\
1\end{array}\right],
\left[\begin{array}{c}
1\\
0\end{array}\right],
\left[\begin{array}{c}
1\\
0\end{array}\right]\right),
\quad
\del_8:=\left(\left[\begin{array}{c}
0\\
1\end{array}\right],
\left[\begin{array}{c}
1\\
0\end{array}\right],
\left[\begin{array}{c}
2\\
0\end{array}\right]\right),
$$
%Non-flat (i.e. non-PRV) solution of {\bf Q4}:
$$
\del_9:=\left(\left[\begin{array}{c}
0\\
1\end{array}\right],
\left[\begin{array}{c}
0\\
1\end{array}\right],
\left[\begin{array}{c}
2\\
0\end{array}\right]\right),
$$
%Not a solution of {\bf Q4}:
$$
\eps_1:=\left(\left[\begin{array}{c}
0\\
1\end{array}\right],
\left[\begin{array}{c}
0\\
1\end{array}\right],
\left[\begin{array}{c}
1\\
0\end{array}\right]\right),
\quad
\eps_2:=\left(\left[\begin{array}{c}
0\\
1\end{array}\right],
\left[\begin{array}{c}
0\\
1\end{array}\right],
\left[\begin{array}{c}
1\\
1\end{array}\right]\right).
$$
\noindent It was observed in \cite{KLM3} that only $\eps_1,  \eps_2$ are not in $Tens(G)$. Moreover,
$$\eps_1+\eps_2=\del_4+\del_9\in Tens(G),\quad
\eps_1+u(\eps_2)=\del_4+u(\del_1)+u(\del_2)\in Tens(G)$$
and for each natural number $n\ge 2$, $n\eps_i$ belongs to $Tens(G)$, $i=1, 2$.

Therefore, if $\si=(\la_1,\la_2,\la_3)$
is a combination of the semigroup generators which is not in $Tens(G)$, then it has the form
\begin{equation}
\label{e1}
\eps_i+ \sum_{g\in S_3} \sum_{j=1}^9 n_{gj} g(\del_j),
\end{equation}
where $n_{gj}\in \Z_+$. By assumption, either $\la_1$ or  $\la_2\notin \Z_+ \varpi_2$.
Therefore at least one of the summands  $g(\del_j)=(\mu_1,\mu_2,\mu_3)$ is
such that $\mu_1$ or $\mu_2$ resp. does not belong to $\Z_+ \varpi_2$ either. Hence Proposition \ref{p0}
would follow from:

\begin{lem}
Suppose that $w(\del_j)=(\mu_1,\mu_2,\mu_3)$ is such that $\mu_1\notin \Z_+\varpi_2$. Then for each $i=1,2$ the sum
$\eps_i+w(\del_j)$ belongs to $Tens(G)$.
\end{lem}
\proof The proof of this lemma is a direct computation with the LiE program. \qed

This concludes the proof of Proposition \ref{p0}. \qed

\begin{rem}
Proposition \ref{p0} implies Corollary \ref{c1}.
\end{rem}

Below we observe that certain combinations of the type (\ref{e1}) not covered by Proposition \ref{p0} nevertheless
belong to $Tens(G)$:

\begin{lem}
The following combinations of the type (\ref{e1}) belong to $Tens(G)$:
$$
\eps_1+\del_5, \quad \eps_1+\del_9,
$$
$$
\eps_2+\del_6,  \quad \eps_2+\del_9,
$$
$$
\eps_1+w(\del_2)+\del_6,  \quad \eps_1+\del_4+\del_6,  \quad \eps_1+u(\del_2)+\del_6,
$$
$$
\eps_2+\del_2+\del_5,  \quad \eps_2+\del_4+\del_5.
$$
\end{lem}
\proof Observe that
$$
\eps_1+\del_5=2\del_9,  \quad \eps_1+\del_9=\del_2+\del_5,
$$
$$
\eps_2+\del_6=\del_9+ u(\del_2),  \quad \eps_2+\del_9=\del_4+\del_5
$$
and hence they are in $Tens(G)$. The combination $\eps_1+u(\del_2)+\del_6$ belongs to $(2L)^3$
and therefore it is in $Tens(G)$. Moreover
$$
\eps_1+w(\del_2)+\del_6=
\eps_2+\del_2+\del_5.
$$
Thus it remains to check 3 last combinations in Lemma, which is done by a direct computation with LiE.
\qed

By combining Proposition \ref{p0} with the above lemma we see that
it remains to analyze combinations of the following types:

\begin{enumerate}
\item $\eps_1+x \del_2+ y w(\del_2) + z\del_4$.

\item $\eps_2+ x\del_2+ yw(\del_2) +z\del_4$.

\item $\eps_1+ n\del_2+ m \del_6$.

\item $\eps_2+ nw(\del_2) + m\del_5$.
\end{enumerate}

We note that the first two types belongs to $\E_1$, the third type belongs to $\E_2$ and the last type belongs to
$\E_3$. Hence we have proved that if $\si\in \P(G)\cap L^3$ does not belong to
$S_3\cdot(\E_1 \cup \E_2\cup \E_3)$ then $\si\in Tens(G)$. It is left to show that $(\E_1 \cup \E_2\cup \E_3) \cap
Tens(G)=\emptyset$. This is done in the following three propositions.

\begin{figure}[tbh]
\centerline{\epsfxsize=2.5in \epsfbox{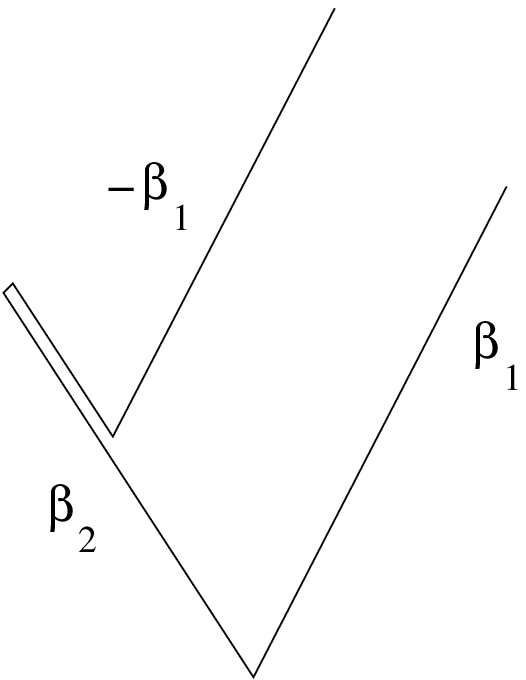}}
\caption{\sl }
\label{f1}
\end{figure}

\begin{prop}\label{p1}
For any $n, m\in \Z_+$, the triple
$$
\si=(\la, \nu, \mu)=\eps_1+ n\del_2+m\del_6=\left(\left[\begin{array}{c}
0\\
1+n+m\end{array}\right],
\left[\begin{array}{c}
0\\
1+n+2m\end{array}\right],
\left[\begin{array}{c}
1+3m\\
0\end{array}\right]\right)
$$
 does not belong to $Tens(G)$. Thus $\E_1\cap Tens(G)=\emptyset$. 
 %Moreover, $\E_1\cap Hecke(G)=\emptyset$. 
\end{prop}
\proof Note that $\la, \nu \in \N\varpi_2, \mu\in \N \varpi_1$. If $\si\in Tens(G)$ then there exists a 
Hecke path $p$ connecting $\la$ to $\nu$ which is entirely contained in $\De$, so that 
$\Length(p)=(1+3m)\varpi_1$. Let $p$ be modeled on a chain $\eta_0\ge \eta_1\ge ... \ge \eta_m$ 
whose elements are in $W(\varpi_1)$.
% determined by the value of the derivative $p'(t_i)$
%where the germ of $p$  at $t_i\in [0,1]$ is geodesic.
Note that since $p(0), p(1)$ are on the wall
$H_2=\R\cdot \varpi_2$, and the image of $p$ is contained in $\De$, it follows that
$$
\eta_0\ne -\be_0, -\be_2, \be_1$$
and
$$
p'(1)=\eta_m\ne \be_2, \be_0.
$$

\begin{figure}[h]
\centerline{\epsfxsize=4.5in \epsfbox{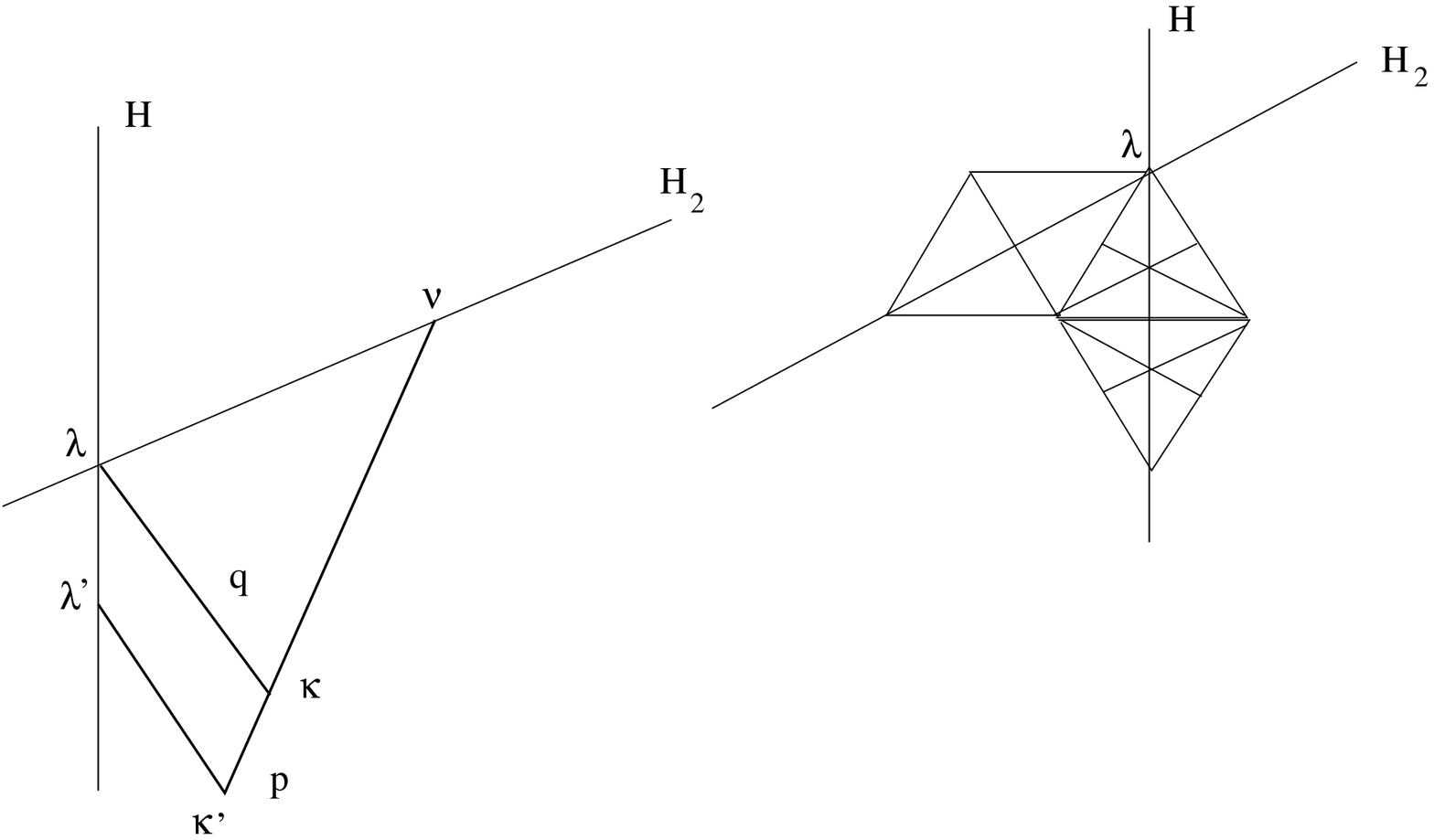}}
\caption{\sl }
\label{f2}
\end{figure}

Therefore, the chain $\eta_0\ge \eta_1\ge ... \ge \eta_m$ is a subchain of
$$
-\be_1\ge -\be_2\ge \be_2\ge \be_1
$$
and the path $p$ has the shape as in Figure \ref{f1}. It is clear that $\eta_m=\be_1$, 
for otherwise the path $p$ cannot connect $\la$ to $\nu$.

We define a {\em canonical path} $q$ which is an $LS(R)$-path connecting $\la$ to $\nu$
and which corresponds to the chain $\be_2\ge \be_1$. Then $\Length(p)=3m\varpi_1$. Let $\kappa$ denote
the break-point of $q$, see Figure \ref{f2}.

\begin{figure}[tbh]
\centerline{\epsfxsize=1.5in \epsfbox{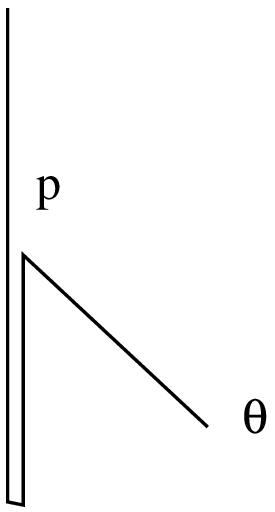}}
\caption{\sl }
\label{f3}
\end{figure}

If $p$ does not contain the subsegment $\ol{\kappa \nu}$, its orthogonal projection to the horizontal
wall $H_1$ will have length strictly less than the length of the orthogonal projection of $\ol{\la \nu}$,
and therefore $p$ cannot connect $\la$ to $\nu$. Thus $p$ contains the subsegment $\ol{\kappa \nu}$,
let $\ol{\kappa'\nu}$ denote the maximal geodesic subsegment in $p$ containing $\ol{\kappa\nu}$. Then
$\kappa'$ has to be a special vertex (since $p$ is a Hecke path). Let $H$ denote
the vertical wall through $\la$.
Let $t$ denote the largest point in $[0,1]$ such that $\la'=p(t)$ belongs to $H$.
Then
$$
\Length(p|[t,1])= \Length(q)=3m\varpi_1.
$$
Suppose that $\kappa'\ne\kappa$.
Since $\kappa'$ is a special vertex, the distance between $\la'$ and $\la$ is at least
the length of $\varpi_2$, i.e. strictly greater than the length of $\varpi_1$. Thus $\kappa'=\kappa$.

Therefore let $\ol{\la'' \kappa}\subset \ol{\la\kappa}$ denote the largest subsegment contained in $p$.
However, if $\la''\ne\la$, we again get a contradiction: The path $p$ is strictly to the right
of the wall $H$ which is absurd. \qed

\begin{prop}
\label{p2}
No triple
$$
\si:=(\la, \nu, \mu)=\left(\left[\begin{array}{c}
0\\
x\end{array}\right],
\left[\begin{array}{c}
0\\
y\end{array}\right],
\left[\begin{array}{c}
1\\
z\end{array}\right]\right)
$$
belongs to $Tens(G)$. Thus $\E_2\cap Tens(G)=\emptyset$.
\end{prop}
\proof We will need the following two lemmas:

Suppose that $\xi\in H_2$ is a special vertex and set $\theta:=\xi+\varpi_1$.

\begin{figure}[tbh]
\centerline{\epsfxsize=4.5in \epsfbox{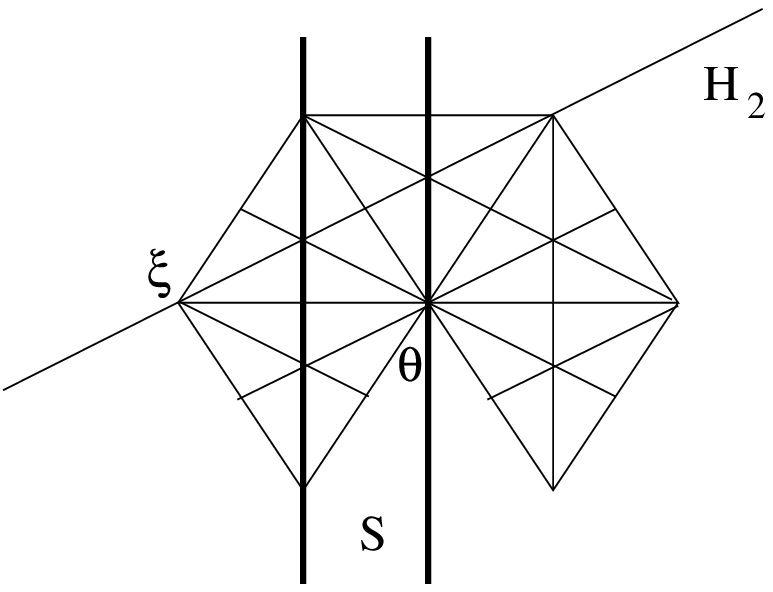}}
\caption{\sl }
\label{f4}
\end{figure}

\begin{lem}\label{pl1}
There are no Hecke paths $p_2: [0,1]\to \De$ so that

\begin{itemize}
\item $p_2(0)=\theta$,

\item $p_2(1)$ is a special vertex in $H_2$.

\item $p_2$ is modeled on a subchain in $-\al_1\ge -\al_2\ge \al_2$.
\end{itemize}
\end{lem}
\proof Under the above assumptions the path $p_2$ has the shape as in Figure \ref{f3}.
Therefore the image of $p_2$ is contained in the vertical strip
$S$ of the width $|\varpi_1|/2$, see Figure \ref{f4}.
However $S\cap H_2$ contains no special vertices. \qed

\medskip
Set $\theta'=\xi+2\varpi_1$.

\begin{figure}[tbh]
\centerline{\epsfxsize=4.5in \epsfbox{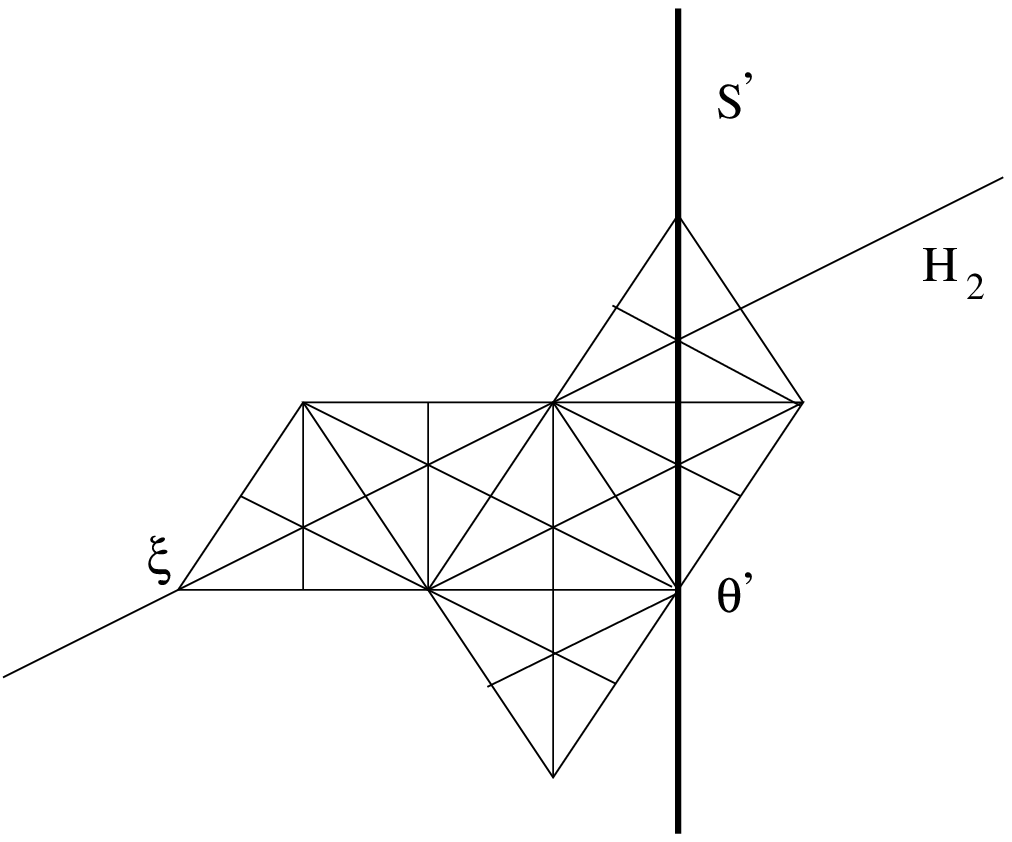}}
\caption{\sl }
\label{f5}
\end{figure}

\begin{lem}\label{pl2}
There are no Hecke paths $p_2: [0,1]\to \De$ so that

\begin{itemize}
\item $p_2(0)=\theta'$,

\item $p_2(1)$ is a special vertex in $H_2$.

\item $p_2$ is modeled on the a subchain in $-\al_2\ge \al_2$.
\end{itemize}
\end{lem}
\proof Under the above assumptions the image of $p_2$ is contained in the vertical wall
$S'$ through the point $\theta'$, see Figure \ref{f5}.
However $S\cap H_2$ is not a special vertex. \qed

\medskip
Suppose that $p$ is an $LS_1(R)$-path in $\De$ with
$$
\length(p)=(\varpi_1, z\varpi_2)
$$
which is contained in $\De$ and connects $\la$ to $\nu$. Let $p=p_1*p_2$ where $p_i$ are
$LS$ paths.

Since $p_2$ is an LS path, if $p_2$ has a break-point on the wall $H_2$, this point has to be special 
(see Observation \ref{ob}). Therefore, we can reduce the discussion to the case when $p_2$ 
does not contain nondegenerate subsegments in $H_2$, which we assume from now on.

\begin{figure}[tbh]
\centerline{\epsfxsize=3.5in \epsfbox{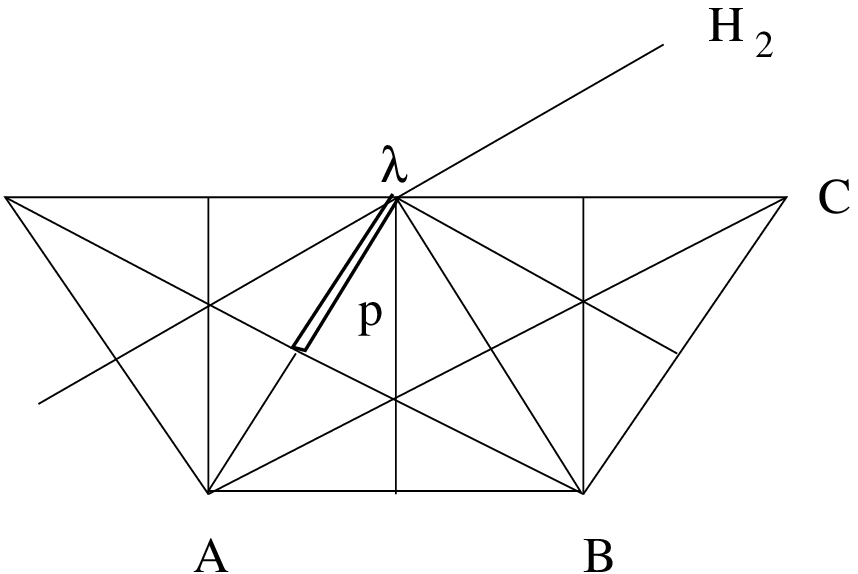}}
\caption{\sl }
\label{f6}
\end{figure}

The path $p_1$ is either a geodesic path connecting $\la$ to $A, B$ or $C$, or $p_1(0)=p_1(1)=\la$
and $p_1$ has a unique break-point which is a point of backtracking, see Figure \ref{f6}.

\begin{lem}
There is no LS path $p_1\subset \De$ connecting $\la$ to itself such that
$$\Length(p_1)=\varpi_1.$$
\end{lem}
\proof Since $p_1(0)=p_1(1)\in H_2$, the path $p_1$ has to be modeled on a  subchain in
$$
-\be_1\ge \be_2\ge \be_1.
$$
Since $\Length(p_1)=\varpi_1$, the path $p_1$ can have only one break-point, hence it is a backtrack.
Thus $p_1$ is modeled on a chain of the form $-\eta\ge \eta$, hence the model chain is $-\be_1\ge \be_1$.
Such a path is Hecke but not an LS path, see Figure \ref{f6}. \qed

\medskip
If $p_1(1)=C$ then $p_1'(1)=\varpi_1$. In this case, $p_2$ is either modeled on the chain
$\al_1\ge \varpi_2$ or $p_2$ is geodesic parallel to the wall $H_2$. In either case,
$p_2(1)\ne H_2$.

\begin{figure}[tbh]
\centerline{\epsfxsize=6in \epsfbox{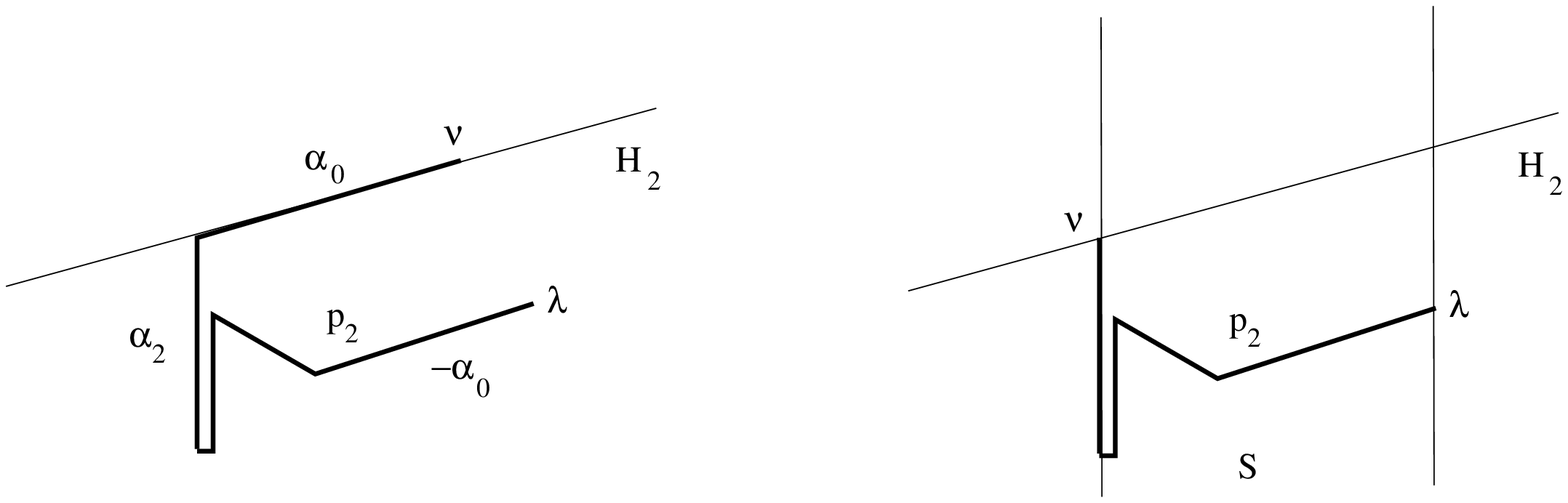}}
\caption{\sl }
\label{f7}
\end{figure}

Thus $p_2(0)$ is either $A$ or $B$. Let $p_2$ be modeled on a chain $\eta_0\ge \eta_1\ge ... \ge \eta_m$
%determined by the values of the derivative $p'(t_i)$, where the germ of $p$ at each
%$t_i\in [0,1]$ is geodesic.
Since $p_2(1)\in H$ and $p_2'(1)\notin \R\varpi_2$, it follows that
$\eta_0\ge ... \ge \eta_m$  is a subchain in
$$
-\varpi_2\ge -\al_1\ge -\al_2 \ge \al_2
$$
and the path $p_2$ has the shape as in the Figure \ref{f7}. In particular, the image of
$p_2$ is contained in the parallel strip $S$ bounded by the vertical walls passing through $p_2(0), p_2(1)$.

\medskip
{\bf Case 1.} $p_2(0)=B$. The the chain condition satisfied by the concatenation $p_1*p_2$ at the
point $B$ implies that the chain $\eta_0\ge ...\ge \eta_m$ is a subchain in $-\al_2\ge \al_2$.
This however contradicts Lemma \ref{pl2}.

\medskip
{\bf Case 2.} $p_2(0)=A$. Then the chain  $\eta_0\ge ...\ge \eta_m$ is a subchain in
$-\al_0\le -\al_1\ge -\al_2\ge \al_2$. If it is a subchain in  $-\al_1\ge -\al_2\ge \al_2$ we get a 
contradiction with Lemma \ref{pl2}. Therefore the initial maximal geodesic segment $\ol{AD}$ in $p_2$ is 
parallel to the wall $H_2$. The assumption that $p_2$ is an LS path then implies that the break-point 
$p_2(t)=D$ is a special vertex. 

\begin{rem}
\label{spec}
There are Hecke paths $p_2: [0,1]\to \De$ for which $p_2(0)=A$, $p_2(1)$ is a special vertex on $H_2$, 
$\Length(p_2)\in \N \varpi_2$, which are modeled on the chain $-\al_0\ge \al_2$. 
However these paths fail to be LS paths, cf. Observation \ref{ob}. 
\end{rem}

Then the restriction $q:=p_2\restr_{[t,1]}$ is a Hecke path in $\De$, so that   
$q(1)=p_2(1)$ a special vertex on $H_2$ and $q$ is modeled on a subchain in $-\al_1\ge -\al_2\ge \al_2$. 
We then obtain a contradiction as above.

This concludes the proof of Proposition \ref{p2}.
\qed
%One can check that $\si$ in Prop \ref{p2} belong to Hecke iff $x, y, z$ satisfy metric triangle inequalities.  
%
%

\begin{prop}\label{p3}
For any $n, m\in \Z_+$, the triple
$$
\si=(\mu, \nu, \la)=\eps_1+ nw(\del_2)+m\del_5=\left(\left[\begin{array}{c}
0\\
1+n+m\end{array}\right],
\left[\begin{array}{c}
0\\
1+m\end{array}\right],
\left[\begin{array}{c}
1+3m\\
1+n\end{array}\right]\right),
$$
 does not belong to $Tens(G)$. Thus $\E_3\cap Tens(G)=\emptyset$.
 %It seems that there are no Hecke triples in this case either. 
 %
\end{prop}
\proof If $\si\in Tens(G)$, then there exits LS path $p\subset \De$  connecting $\la$ to $\nu$ so that
$\Length(p)=\mu$. The path $p$ is modeled on a subchain in the chain
$$
-\varpi_2\ge -\al_1\ge -\al_2\ge \al_2\ge \om_2.
$$
Therefore the general shape of $p$ is as in Figure \ref{f7}. Therefore $p$ lies to the right of the
vertical wall passing through its last break-point $p(t)$. Similarly to the proof of Proposition \ref{p2} 
it suffices to consider the case when $p$ contains no nondegenerate subsegments of the wall $H_2$ 
(see Figure \ref{f7}), i.e. 
$p$ is modeled on a subchain in 
$$
-\varpi_2\ge -\al_1\ge -\al_2\ge \al_2.
$$
Let $H$ denote the vertical wall through $\nu$. Set $\xi:= \nu-\frac{2}{3}\varpi_2$,
$\zeta:= \xi+\varpi_2$ and let $H'$ denote the vertical wall through $\zeta$.
Let $\theta\in H$ denote the point so that the triangle $[\xi, \theta, \nu]$
is equilateral. Thus $\nu$ belongs to the interior of the segment $\ol{\xi \zeta}$.

\begin{figure}[tbh]
\centerline{\epsfxsize=5in \epsfbox{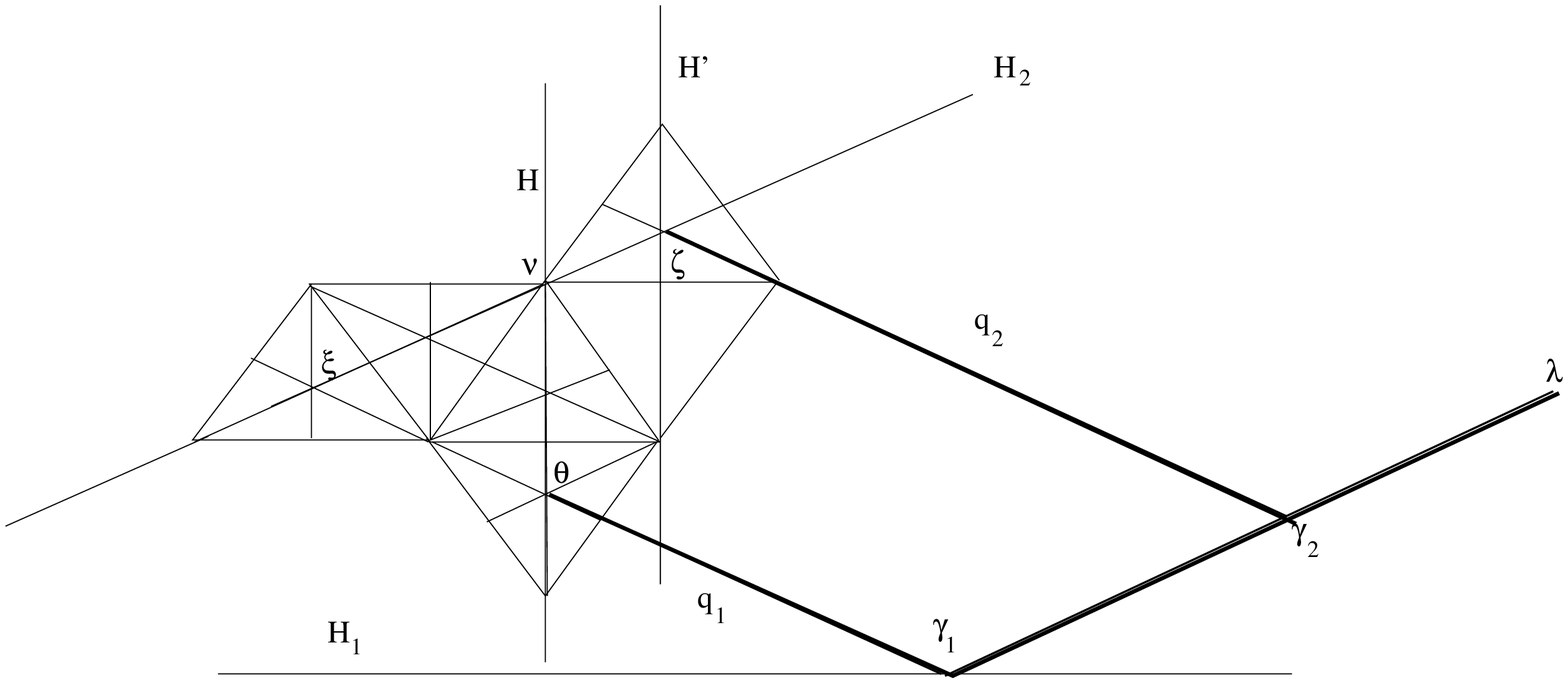}}
\caption{\sl }
\label{f8}
\end{figure}

Define paths $q_1, q_2$ from $\la$ to $\theta$ and $\zeta$ as in Figure \ref{f8}. These
paths have break-points at $\ga_1=(1+3m)\varpi_1$ and $\ga_2=(1+3m)\varpi_1+\varpi_2$ respectively.
Then
$$
\Length(q_1)=(n+m+\frac{2}{3})\varpi_2, \Length(q_2)=(n+m+\frac{2}{3})\varpi_2.
$$
Note that $q_2$ lies entirely to the right of the wall $H'$.

If $p$ does not contain the subsegment $\ol{\la \ga_2}$, it also lies to the right of the wall
$H'$. Such a path cannot connect $\la$ to $\nu$.

Therefore, $p$ contains the segment $\ol{\la \ga_1}$. The same argument shows that $p$ cannot have the first
break at $\ga_2$. Moreover, unless $p$ is modeled on a chain of the form
$-\varpi_2 \ge -\al_1 \ge ...$, it lies to the right of the vertical wall through $\ga_1$.
Thus, since $p$ is an LS path, its first break-point has to be a special vertex.

%\begin{rem}
%Analogously to Remark \ref{spec}, in the case when $p$ is a Hecke path 
%modeled on the chain $-\al_0\ge \al_2$, the break-point of $p$ can be nonspecial. 
%\end{rem}

 The interior of
the segment $\ol{\ga_1 \ga_2}$ contains no special vertices, hence the first break of $p$ occurs
at $\ga_1$ and $p$ contains a subsegment of $\ol{\ga_1 \theta}$.

If $p$ makes a break after $\ga_1$ but before reaching
the wall $H$, it lies strictly to the right of $H$ which is again impossible. Thus $p$ contains
the subpath $q_1$; let $t\in [0,1]$ be such that $p(t)=\theta$. Then
$$
\Length(p|[t,1])= \mu-\Length(q_1)=\frac{1}{3}\varpi_2.
$$
However, the distance from $\theta$ to $\nu$ equals $\frac{2}{3}|\varpi_2|>\frac{1}{3}|\varpi_2|$.
Thus $p(1)\ne \nu$. Contradiction. \qed

This concludes the proof of Theorem \ref{g2}.

%\newpage


\begin{thebibliography}{BaBE}
\addcontentsline{toc}{section}{Bibliography}


\bibitem{Ballmann}
W.\ Ballmann, ``Lectures on spaces of nonpositive curvature.
With an appendix by Misha Brin.'' DMV Seminar, vol. 25.
Birkhauser Verlag, Basel, 1995.

\bibitem{BS}
A.\ Berenstein and R.\ Sjamaar,
{\em Coadjoint orbits, moment polytopes,
and the Hilbert-Mumford criterion},  Journ. Amer. Math. Soc., vol. {\bf 13} (2000), no. 2, p. 433--466.



\bibitem{BZ} A. Berenstein and  A. Zelevinsky,
{\em Tensor product multiplicities, canonical bases and totally
positive varieties,}  Invent. Math.  143  (2001),  no. 1, p. 77--128.


\bibitem{Bourbaki}
N. Bourbaki, ``Lie groups and Lie algebras", Chap. 4, 5, 6. Springer Verlag, 2002.

\bibitem{Brown}
K. Brown,``Buildings'', Springer-Verlag, New York, 1989.


\bibitem{BT}
F.\ Bruhat\ and\ J. Tits, {\em Groupes reductifs sur un corps
local,} Inst. Hautes \'Etudes Sci. Publ. Math., No. 41, (1972), p.
5--251.


\bibitem{C}
R. Cluckers, {\em Presburger sets and $p$-minimal fields}, J.
Symbolic Logic, 68 (2003), no. 1, 153--162.


\bibitem{FultonHarris}
W.\ Fulton and J.\ Harris, ``Representation theory. A first
course,'' Graduate Texts in Math., vol.  {\bf 129}, Springer
Verlag, 1991.

\bibitem{H}
J. E. Humphreys,
``Introduction to Lie algebras and representation theory.''
Graduate Texts in Mathematics, 9. Springer-Verlag, New York-Berlin, 1978.



\bibitem{KLM1}
M. Kapovich, B.\ Leeb and J.\ J.\ Millson, {\em Convex functions on
symmetric spaces, side lengths of polygons and the stability
inequalities for weighted configurations at infinity}, Preprint,
June 2004.

\bibitem{KLM2}
M.\ Kapovich, B.\ Leeb and J. J. Millson, {\em Polygons in buildings
and their refined side-lengths}, Preprint, June 2004,
arXiv:math.MG/0406305.


\bibitem{KLM3}
M. Kapovich,  B. Leeb and J. J. Millson, {\em Polygons in
symmetric spaces and buildings with applications to
  algebra},  Preprint, 2004, arXiv:math.RT/0210256, To appear in Memoirs of AMS.



\bibitem{KM}
M. Kapovich and  J. J. Millson, {\em A path model for geodesics in
Euclidean buildings and its applications to  representation
theory}, Preprint, 2004, arXiv:math.RT/0411182, Submitted to Math. Publ. of IHES.

\bibitem{KleinerLeeb}
B.\ Kleiner and B.\ Leeb, {\em Rigidity of quasi-isometries for
symmetric spaces and Euclidean buildings}, Publ. Math. IHES, vol.
{\bf 86} (1997), p. 115--197.


\bibitem{KT}
A. Knutson\ and\ T. Tao, {\em The honeycomb model of ${\rm GL}_n(\C)$ tensor products. I.
Proof of the saturation conjecture,} J. Amer. Math. Soc., vol. {\bf 12} (1999), no.~4, p. 1055--1090.



\bibitem{KuLM}
S.\ Kumar, B.\ Leeb and J.\ J.\ Millson, {\em The generalized
triangle inequalities for rank $3$ symmetric space of noncompact
type}, In: ``Explorations in complex and Riemannian geometry
(Papers dedicated to Robert Greene),'' Contemporary Math., Vol
332, 2003, p. 171--195.



\bibitem{Laskowski}
M.\ C.\ Laskowski, {\em An application of Kochen's Theorem,} J.
Symbolic Logic, {\bf 68} (2003), no. 4, 1181--1188.



\bibitem{Littelmann2}
P.\ Littelmann, {\em Paths and root operators in representation
theory}, Annals of Math. (2)  142 (1995) no. 3, p. 499--525.


\bibitem{Ronan}
M. Ronan,``Lectures on buildings'', Perspectives in Mathematics, 7.
Academic Press, Inc., 1989.


\bibitem{Rousseau}
G. Rousseau, {\em Euclidean buildings}, Lectures at Ecole d'ete de
Mathematiques ``Nonpositively curved geometries, discrete groups and
rigidities'', Institut Fourier, Grenoble, 2004.


\bigskip
\noindent Michael Kapovich: \newline Department of Mathematics,
\newline University of California, \newline Davis, CA 95616, USA
\newline
 kapovich$@$math.ucdavis.edu


\smallskip
\noindent John J. Millson: \newline Department of Mathematics,
\newline University of Maryland, \newline College Park, MD 20742,
USA
\newline
 jjm$@$math.umd.edu


\end{thebibliography}
\end{document}